\documentclass[12pt]{article}

\usepackage{amsmath}
\usepackage{amssymb}
\usepackage{vatola}
\usepackage{CJK}
\usepackage{amsmath}
\usepackage{amsthm}
\usepackage{mathrsfs}
\usepackage{graphicx}
\usepackage{amssymb}
\usepackage{latexsym}
\allowdisplaybreaks[4] 
\def\q{\quad}
\def\qq{\qquad}
\def\t{\hbox}
\def\e{\equiv}
\def\f{\frac}
\def\b{\binom}
\def\ex{\text{ex}}

\let\pro=\proclaim
\let\endpro=\endproclaim

\def\qtq#1{\q\t{#1}\q}
\theoremstyle{plain}

\theoremstyle{definition}

\begin{document}
\begin{CJK*}{GBK}{song}
 \centerline {\bf Tur\'an's problem for trees $T_n$ with maximal degree
 $n-4$}
\par\q\newline \centerline{Zhi-Hong Sun$^1$ and Yin-Yin Tu$^2$} \par\q\newline \centerline{$\ ^1$School of Mathematical
Sciences, Huaiyin Normal University,} \centerline{Huaian, Jiangsu
223001, P.R. China} \centerline{Email: zhihongsun@yahoo.com}
\centerline{Homepage: http://www.hytc.edu.cn/xsjl/szh}

\centerline{$\ ^2$School of Mathematical Sciences, Huaiyin Normal
University,} \centerline {Huaian, Jiangsu 223001, PR China}
\centerline{Email: yinyintu@126.com}

\abstract{For $n\ge 6$ let $V=\{v_0,v_1,\ldots,v_{n-1}\}$,
$E_1=\{v_0v_1,\ldots,v_0v_{n-4},v_1v_{n-3},v_1v_{n-2}$,
$v_1v_{n-1}\}$,
$E_2=\{v_0v_1,\ldots,v_0v_{n-4},v_1v_{n-3},v_1v_{n-2},v_2v_{n-1}\}$,
$E_3=\{v_0v_1,\ldots,v_0v_{n-4}$,
$v_1v_{n-3},v_2v_{n-2},v_3v_{n-1}\}$, $T_n^3=(V,E_1),\
T_n^{''}=(V,E_2)$ and $T_n^{'''} =(V,E_3).$
  In this paper, for $p\ge n\ge 15$ we obtain explicit
formulas for $ex(p;T_n^3)$,  $ex(p;T_n^{''})$ and $ex(p;T_n^{'''})$,
where $ex(p;L)$ denotes the maximal number of edges in a graph of
order $p$ not containing $L$ as a subgraph.
\par\q
\newline MSC: Primary 05C35, Secondary 05C05
\newline Keywords: tree, Tur\'an problem} \endabstract
 \footnotetext[1] {The first author is supported by the
National Natural Science Foundation of China (grant No. 11371163).}

\section*{1. Introduction}
\par\q
\par In this paper, all graphs are simple graphs. For a graph $G=(V(G),E(G))$
let $e(G)=|E(G)|$ be the number of edges in $G$ and let $\Delta(G)$
be the maximal degree of $G$.
  For a forbidden graph $L$, let $ex(p;L)$
denote the maximal number of edges in a graph of order $p$ not
containing $L$ as a subgrph. The corresponding Tur\'an's problem is
to evaluate $ex(p;L)$.
\par Let $\Bbb N$ be the set of
positive integers, and let $p,n\in\Bbb N$ with $p\ge n\ge 3$. For a
given tree $T_n$ on $n$ vertices,
 it is difficult to determine the value of $ex(p;T_n)$. The famous Erd\"os-S\'os
 conjecture asserts that
$ex(p;T_n)\le \f{(n-2)p}2$ for every tree $T_n$ on $n$ vertices. For
the progress on the Erd\"os-S\'os
 conjecture, see for example [2,5].
 Write
$p=k(n-1)+r$, where $k\in\Bbb N$ and $r\in\{0,1,\ldots,n-2\}$. Let
$P_n$ be the path on $n$ vertices. In [1] Faudree and Schelp showed
that
$$ex(p;P_n)=k\binom {n-1}2+\binom r2=\f{(n-2)p-r(n-1-r)}2.\tag 1.1$$
 Let $K_{1,n-1}$ denote the unique tree on $n$
vertices with $\Delta(K_{1,n-1})=n-1$, and for $n\ge 4$ let
  $T_n'$ denote the unique tree on $n$ vertices
with $\Delta(T_n')=n-2$. In [3] the first author and Lin-Lin Wang
determined $ex(p;K_{1,n-1})$ and $ex(p;T_n')$. In [3,4] the first
author and his coauthors also determined $ex(p;T_n)$ for trees $T_n$
with $n$ vertices and $\Delta(T_n)=n-3$.
\par For $n\ge 6$ let
$$\align &V=\{v_0,v_1,\ldots,v_{n-1}\},\q
E_1=\{v_0v_1,\ldots,v_0v_{n-4},v_1v_{n-3},v_1v_{n-2},v_1v_{n-1}\},
\\&E_2=\{v_0v_1,\ldots,v_0v_{n-4},v_1v_{n-3},v_1v_{n-2},v_2v_{n-1}\},
\\&E_3=\{v_0v_1,\ldots,v_0v_{n-4},v_1v_{n-3},v_2v_{n-2},v_3v_{n-1}\}.
\endalign$$ Suppose $T_n^3=(V,E_1),\ T_n^{''}=(V,E_2)$ and $T_n^{'''}
=(V,E_3).$
  In this paper, for $p\ge n\ge 15$ we obtain explicit
formulas for $ex(p;T_n^3)$,  $ex(p;T_n^{''})$ and $ex(p;T_n^{'''})$,
see Theorems 3.1, 5.1 and 4.1-4.5.
\par
In addition to the above notation, throughout this paper we also use
the following notation: $[x]\f{\q}{\q}$the greatest integer not
exceeding $x$, $d(v)\f{\q}{\q}$the degree of the vertex $v$ in a
graph,  $d(u,v)\f{\q}{\q}$the distance between the two vertices $u$
and $v$ in a graph,  $K_n\f{\q}{\q}$the complete graph on $n$
vertices, $K_{m,n}\f{\q}{\q}$the complete bipartite graph with $m$
and $n$ vertices in the bipartition,
 $\overline{G}\f{\q}{\q}$the complement of $G$,  $G[V_1]\f{\q}{\q}$the subgraph
of $G$ induced by vertices in the set $V_1$, $G-V_1\f{\q}{\q}$the
subgraph of $G$ obtained by deleting vertices in $V_1$ and all edges
incident with them, $\Gamma(v)\f{\q}{\q}$the set of vertices
adjacent to the vertex $v$, $\Gamma_2(v)\f{\q}{\q}$the set of those
vertices $u$ such that $d(u,v)=2$, $e(V_1V_1')\f{\q}{\q}$the number
of edges with one endpoint in $V_1$ and another endpoint in $V_1'$.

\section*{2. Basic lemmas}\par\q
\pro{Lemma 2.1} Let $p,n \in \Bbb N$ with $p\geq n\geq 10$. Let
$T_n$ be a tree with $n$ vertices
 and $\Delta(T_n)=n-4$, and let $G\in
Ex(p;T_n)$.  Then $\Delta(G)\ge n-5$.\endpro
 Proof. By [3, Theorem 2.1], $ex(p;K_{1,n-4})=[\f{(n-5)p}2]$.
Since a graph does not contain $K_{1,n-4}$ as a subgraph implies
that the graph does not contain any copies of $T_n$, we have
$$e(G)=ex(p;T_n)\geq ex(p;K_{1,n-4})=\big[\f{(n-5)p}2\big].
\tag 2.1$$
If $\Delta(G)\leq n-6$, using Euler's theorem we see that
$e(G)=\f12\sum_{v\in V(G)}d(v)\le \f{(n-6)p}2$. Hence $\f{(n-5)p-1}2
\leq [\f{(n-5)p}2]\leq e(G)\leq \f{(n-6)p}2$. This is impossible.
Thus $\Delta(G)\geq n-5$.
 \pro{Lemma 2.2} Let $p,n \in
\Bbb N$ with $p\geq n\geq 10$. Let $T_n$ be a tree with $n$ vertices
and
 $G\in Ex(p;T_n)$.  Suppose $V_1\subset V(G)$ and $|V_1|=m+1\ge
n-3$. Then
 $e(G)-e(G-V_1)>3m$.\endpro
Proof. We first assume $m\ge n-2$. Suppose $m+1=k(n-1)+r$ with $k\in
\Bbb N$ and $r\in \{0,1,\cdots,n-2\}$.  Then clearly $kK_{n-1}\cup
K_r$ does not contain any copies of $T_n$ and
$$\aligned  e(kK_{n-1}\cup
K_r)&= \f{k(n-1)(n-2)}2+\f{r(r-1)}2 =\f{(n-2)(m+1)-r(n-1-r)}2
\\&\geq \f{(n-2)(m+1)}2-\f{(n-1)^2}8.\endaligned$$
 As $n\ge 10$ and $m+1\ge n-1$ we have
 $(n-8)(m+1)\ge (n-8)(n-1)>\f{(n-1)^2}4-6$ and so
$$ex(m+1;T_n)\ge e(kK_{n-1}\cup
K_r)\geq \f{(n-2)(m+1)}2-\f{(n-1)^2}8>3m.\tag 2.2$$ If
$e(G)-e(G-V_1)\le 3m$, then
$$e(G)< e(G-V_1)+e(kK_{n-1}\cup K_r)=e((G-V_1)\cup kK_{n-1}\cup
K_r).$$ This contradicts the assumption $G\in Ex(p;T_n)$. If $m=n-3$
or $n-4$, then $e(K_{m+1}) =\f{m(m+1)}2>3m$. If $e(G)-e(G-V_1)\le
3m$, then $e(G)< e(G-V_1)+e(K_{m+1})=e((G-V_1)\cup K_{m+1}),$ which
contradicts the assumption $G\in Ex(p;T_n)$. Hence
$e(G)-e(G-V_1)>3m$ as claimed.

\pro{Lemma 2.3} Let $p,n \in \Bbb N, p\geq n\geq 10$, $p=k(n-1)+r$,
 $k\in\Bbb N$ and $r\in\{0,1,\ldots,n-2\}$.
 Let $T_n$ be a tree with n vertices,
 $\Delta(T_n)=n-4$ and $G\in
Ex(p;T_n)$. If $G$ is connected and $\Delta(G)\le n-4$, then $p\le
\t{min}\{\f {3(n-1+r)+((-1)^n+(-1)^r)/2}2,\f{r(n-1-r)}2\}$ and so
$p\le 2n-7$.\endpro Proof. Suppose that $G$ is connected and
$\Delta(G)\le n-4$. By [3, Theorem 2.1],
$ex(n-1+r;K_{1,n-4})=[\f{(n-1+r)(n-5)}2]$. Let $G_0\in
Ex(n-1+r;K_{1,n-4})$. Then $G_0$ does not contain $T_n$ and so
$(k-1)K_{n-1}\cup G_0$ does not contain $T_n$ as a subgraph. Thus,
$$e((k-1)K_{n-1}\cup G_0)\le ex(p;T_n)=e(G)\le\Big[\f{(n-4)p}2\Big]
=\f{(n-4)p-(1-(-1)^{nr})/2}2.$$ On the other hand,
$$\aligned &e((k-1)K_{n-1}\cup G_0)\\&=(k-1)\b{n-1}2+\Big[\f{(n-1+r)(n-5)}2\Big]
\\&= (k-1)\b{n-1}2+\f{(n-1+r)(n-5)-(1-(-1)^{(n-1)(r-1)})/2}2
\\&=\f{(n-4)p-(1-(-1)^{nr})/2}2+p-\f{3(n-1+r)+((-1)^{nr}-(-1)^{(n-1)(r-1)})/2}2
\\&=\f{(n-4)p-(1-(-1)^{nr})/2}2+p-\f{3(n-1+r)+((-1)^n+(-1)^r)/2}2
.\endaligned$$ Thus, $p\le \f {3(n-1+r)+((-1)^n+(-1)^r)/2}2$. We
also have
$$\f{(n-2)p-r(n-1-r)}2=e((k-1)K_{n-1}\cup K_r)\le e(G)\le
\f{(n-4)p}2$$ and so $p\le \f{r(n-1-r)}2$. Hence $p\le \t{min}\{\f
{3(n-1+r)+((-1)^n+(-1)^r)/2}2,\f{r(n-1-r)}2\}$.
\par  As $p\ge n$,
we see that $r\not\in\{0,1,2,n-3,n-2\}$ and so $p\le
\f{3(n-1+n-4)+1}2=3n-7$. If $p\ge 2(n-1)$, then $p=2(n-1)+r$ with
$0\le r \le n-5$. As $2(n-1)+r>\f{3(n-1+r)+1}2$, we get a
contradiction. Hence $p<2n-2$. Now we have $p=n-1+r$ with $0\le r\le
n-4$. As $\Delta(G)\le n-4$ we have $$ex(2n-5;T_n)=e(G)\le
\f{(n-4)(2n-5)}2=n^2-\f{13}2n+10<n^2-6n+11=e(K_{n-1}\cup K_{n-4}),$$
which is a contradiction. Hence $p\le 2n-6$. As
$$ex(2n-6;T_n)=e(G)\le
\f{(n-4)(2n-6)}2=n^2-7n+12<n^2-7n+16=e(K_{n-1}\cup K_{n-5}),$$ we
get $p\le 2n-7$. This proves the lemma.

\pro{Lemma 2.4 ([4, Lemma 2.4])} Let $n,n_1,n_2\in\Bbb N$ with
$n_1<n-1$ and $n_2<n-1.$
\par $(\t{\rm{i}})$ If $n_1+n_2<n,$ then
$\b{n_1}2+\b{n_2}2<\b{n_1+n_2}2$.
\par $(\t{\rm{ii}})$ If $n_1+n_2\ge n,$ then
$\b{n_1}2+\b{n_2}2<\b{n-1}2+\b{n_1+n_2-n+1}2.$ \endpro

\pro{Lemma 2.5} Let $n\in \Bbb N$ with $n\ge 10$, and let $T_n$ be a
tree with n vertices and
 $\Delta(T_n)=n-4$. Suppose that for any
 positive integer $m\ge n$ and connected
 graph $H\in \t{Ex}(m;T_n)$ we have $\Delta(H)\le n-4$.
 Let $p\in \Bbb N$, $p\geq n$,
$p=k(n-1)+r$, where
 $k\in\Bbb N$ and $r\in\{0,1,\ldots,n-2\}$.
Assume that $G\in Ex(p;T_n)$, $G$ is not connected, $G_1,\ldots,
G_s$ are distinct components of $G$, $|V(G_i)|=p_i$
$(i=1,2,\ldots,s)$ and $p_1\le p_2\le\cdots \le p_s$. Then $p_1\le
p_2=\cdots =p_{s-1}=n-1\le p_s\le 2n-7$. If $p_1<n-1$ and $p_s\ge
n$, then $p_1\le n-7$ and $p_1(n-3-p_1)\le p_1+p_s+1\le 2n-7$.
\endpro
Proof.  Suppose that $p_s\ge p_{s-1}\ge n$. Then clearly $G_{s-1}\in
\t{Ex}(p_{s-1};T_n)$, $G_s\in \t{Ex}(p_s;T_n)$ and $G_{s-1}\cup
G_s\in Ex(p_{s-1}+p_s;T_n)$. By the assumption, $\Delta(G_{s-1})\le
n-4$ and $\Delta(G_s)\le n-4$. Hence,
$$\ex(p_{s-1}+p_s;T_n)=e(G_{s-1}\cup G_s)
\le \f{(n-4)(p_{s-1}+p_s)}2.$$ If $p_{s-1}+p_s<3(n-1)-1$ and $G_0\in
Ex(p_{s-1}+p_s-(n-1);K_{1,n-4})$, then $K_{n-1}$ does not contain
$T_n$ and
$$\aligned e(K_{n-1}\cup G_0)&=\b{n-1}2+\Big[\f{(n-5)(p_{s-1}+p_s-(n-1))}2\Big]
\\&\ge
\f{(n-1)(n-2)}2+\f{(n-5)(p_{s-1}+p_s-(n-1))-1}2\\&=\f{(n-5)(p_{s-1}
+p_s)+3(n-1)-1}2
\\&>\f{(n-4)(p_{s-1}+p_s)}2\ge ex(p_{s-1}+p_s;T_n).
\endaligned$$
This contradicts the fact that $K_{n-1}\cup G_0$ does not contain
$T_n$. Hence $p_{s-1}+p_s\ge 3(n-1)-1$.  By Lemma 2.3, $p_{s-1}\le
2n-7$ and $p_s\le 2n-7$. Thus,
$$3(n-1)-1\le p_{s-1}+p_s\le 2(2n-7)< 6(n-1)-1.$$   Suppose $G_0\in
Ex(p_{s-1}+p_s-2(n-1);K_{1,n-4})$. Then $G_0$ does not contain $T_n$
and $e(G_0)=[\f{(n-5)(p_{s-1}+p_s-2(n-1))}2]$. Thus,
$$\aligned &ex(p_{s-1}+p_s;T_n)\\&\ge e(2K_{n-1}\cup G_0)=(n-1)(n-2)
+\Big[\f{(n-5)(p_{s-1}+p_s-2(n-1))}2\Big]\\&=3(n-1)+\Big[\f{(n-5)(p_i+p_j)}2
\Big]\ge 3(n-1)+\f{(n-5)(p_{s-1}+p_s)-1}2
\\&=\f{(n-4)(p_{s-1}+p_s)+6(n-1)-1-(p_{s-1}+p_s)}2>\f{(n-4)(p_{s-1}+p_s)}2.
\endaligned$$
This is a contradiction.
\par By the above,  $p_1\le p_2\le \cdots \le
p_{s-1}\le n-1$.  We claim that $p_2\ge n-1$. Otherwise, $p_1\le
p_2<n-1$ and $G_1\cup G_2\cong K_{p_1}\cup K_{p_2}.$ If $p_1+p_2<n,$
by Lemma 2.4(i) we have
 $$e(G_1\cup G_2)=e(K_{p_1}\cup
 K_{p_2})=\b{p_1}2+\b{p_2}2<\b{p_1+p_2}2=e(K_{p_1+p_2}).$$
Since $K_{p_1+p_2}$ does not contain $T_n$ and $G_1\cup G_2\in
Ex(p_1+p_2;T_n)$ we get a contradiction. Hence $p_1+p_2\ge n$. Using
Lemma 2.4(ii) we see that
 $$\aligned
 e(G_1\cup G_2)&=e(K_{p_1}\cup K_{p_2})=\b{p_1}2+\b{p_2}2
 \\&<\b{n-1}2+\b{p_1+p_2-n+1}2=e(K_{n-1}\cup K_{p_1+p_2-n+1}).
 \endaligned$$
Since $p_1\le p_2<n-1,$ we have $p_1+p_2-n+1<n-1$. Hence
$K_{n-1}\cup K_{p_1+p_2-n+1}$ does not contain $T_n$. As $G_1\cup
G_2$ is an extremal graph without $T_n$, this is a contradiction.
Thus, $p_2\ge n-1$. Hence $p_1\le n-1=p_2=\cdots=p_{s-1}\le p_s\le
2n-7$.
\par Assume that $p_s\ge n$ and $p_1<n-1$.
If $p_1+p_s\ge 2n-5$, setting $G_0\in Ex(p_1+p_s-(n-1);K_{1,n-4})$
we find that  $G_0$ does not contain $T_n$ and
$e(G_0)=[\f{(n-5)(p_1+p_s-(n-1)}2]$. Thus,
$$\aligned &e(K_{n-1}\cup G_0)\\&=\b{n-1}2+\Big[\f{(n-5)(p_1+p_s-(n-1))}2\Big]
=\Big[\f{(n-5)(p_1+p_s)+3(n-1)}2\Big]
\\&\ge \f{(n-5)(p_1+p_s)+3(n-1)-1}2
\\&=\f{(n-4)p_s}2+\b{p_1}2+\f{3(n-1)-1+p_1(n-4-p_1)-p_s}2
\\&>\f {(n-4)p_s}2 +\b{p_1}2\ge e(G_1\cup G_s).\endaligned$$
This contradicts the fact $G_1\cup G_s\in Ex(p_1+p_s;T_n)$. Hence
$p_1+p_s\le 2n-6$. If $p_1\ge n-5$, then $p_s\le 2n-6-p_1\le
2n-6-(n-5)=n-1$. This contradicts the assumption $p\ge n$. Hence
$p_1\le n-6$. If $p_1=n-6$, then $p_s=n$. As
$$\aligned e(K_{n-1}\cup K_{n-5})&=\f{(n-1)(n-2)}2+\f{(n-5)(n-6)}2
>\f{(n-6)(n-7)}2+\f {n(n-4)}2\\&\ge e(G_1\cup G_s)=ex(p_1+p_s;T_n)
,\endaligned$$ we get a contradiction. Hence $p_1\le n-7$. We claim
that $p_s\ge p_1(n-4-p_1)-1$. Otherwise, for  $G_0\in
Ex(p_1+p_s;K_{1,n-4})$ we have
$$\aligned e(G_0)&=\Big[\f{(n-5)(p_1+p_s)}2\Big]\ge
\f{(n-5)(p_1+p_s)-1}2
\\&>\f{(n-4)p_s}2+\f{p_1(p_1-1)}2\ge e(G_1\cup
G_s)=ex(p_1+p_s;T_n),\endaligned$$ which is a contradiction. Hence
the claim is true. As $p_1+p_s\le 2n-6$, we get $p_1(n-4-p_1)-1\le
p_s\le 2n-6-p_1$ and so $p_1(n-3-p_1)\le p_1+p_s+1$. By Lemma 2.1,
$\Delta(G_s)\le n-4$. Thus,
$$ex(p_1+p_s;T_n)=e(G_1\cup G_s)\le \f{p_1(p_1-1)}2+\f{(n-4)p_s}2
=\f{(p_1+p_s)(n-4)-p_1(n-3-p_1)}2.$$ On the other hand,
$$\aligned ex(p_1+p_s;T_n)&\ge e(K_{n-1}\cup K_{p_1+p_s-(n-1)})
\\&=\f{(n-1)(n-2)+(p_1+p_s-(n-1))(p_1+p_s-n)}2
\\&=(n-1)^2-\f{(p_1+p_s)(3n-5-p_1-p_s)+(p_1+p_s)(n-4)}2.
\endaligned$$
Hence $-p_1(n-3-p_1)\ge 2(n-1)^2-(p_1+p_s)(3n-5-p_1-p_s)$ and so
$$\aligned &(p_1+p_s)(3n-5-p_1-p_s)\\&\ge 2(n-1)^2+p_1(n-3-p_1)
\ge 2(n-1)^2+n-4=2n^2-3n-2.\endaligned$$ As
$(2n-6)(3n-5-(2n-6))=2n^2-4n-6<2n^2-3n-2$ and
$(2n-7)(3n-5-(2n-7))=2n^2-3n-14<2n^2-3n-2$, we get
$p_1+p_s\not=2n-6,2n-7$ and so $p_1+p_s\le 2n-8$. Thus
$p_1(n-3-p_1)\le p_1+p_s+1\le 2n-7$. This completes the proof.

\pro{Lemma 2.6} Let $n\in \Bbb N$ with $n\ge 10$, and let $T_n$ be a
tree with n vertices and
 $\Delta(T_n)=n-4$. Suppose that for any positive integer
 $m\ge n$ and connected
 graph $H\in \t{Ex}(m;T_n)$ we have $\Delta(H)\le n-4$.
 Let $p\in \Bbb N$ with $p\ge 2n-6$. Then
$$ex(p;T_n)=\f{(n-1)(n-2)}2+ex(p-(n-1);T_n).$$
\endpro
Proof. Let $G\in Ex(p;T_n)$. As $p\ge 2n-6>2n-7$, we see that $G$ is
not connected by Lemma 2.3. Suppose that $G_1,\cdots,G_s$ are all
distinct components of G with $|V(G_i)|=p_i$ and $p_1\leq p_2\leq
\cdots \leq p_s$. Then clearly $G_i\in Ex(p_i;T_n)$ for
$i=1,2,\cdots,s$. By Lemma 2.5, $p_1\le p_2=\cdots =p_{s-1}=n-1\le
p_s\le 2n-7$. If $p_i=n-1$ for some $i\in\{1,2,\ldots,s\}$, then
clearly the result holds. If $p_i\not=n-1$ for all $i=1,2,\ldots$,
then $s=2$, $p_1<n-1<n\le p_2$. By Lemma 2.5, $p=p_1+p_2\le 2n-8$,
which contradicts the assumption $p\ge 2n-6$. Hence the theorem is
proved.

\pro{Lemma 2.7} Let $n\in \Bbb N$ with $n\ge 10$, and let $T_n$ be a
tree with n vertices and
 $\Delta(T_n)=n-4$. Suppose that for any positive integer
 $m\ge n$ and connected
 graph $H\in \t{Ex}(m;T_n)$ we have $\Delta(H)\le n-4$.
 Assume $p,k\in \Bbb N$, $p=k(n-1)+r$, $k\ge 2$ and
$r\in\{0,1,\ldots,n-2\}$. Then
$$\ex(p;T_n)=\f{(n-2)(p-(n-1+r))}2+\ex(n-1+r;T_n).$$
\endpro
Proof. By Lemma 2.6,
$$\align &\t{ex}(p;T_n)\\&=\sum_{s=2}^k\big(\ex(s(n-1)+r;T_n)-\ex((s-1)(n-1)+r;T_n)\big)
+\ex(n-1+r;T_n)
\\&=(k-1)\b {n-1}2+\ex(n-1+r;T_n).\endalign$$
Since $(k-1)(n-1)=p-(n-1+r)$ we deduce the result.

\pro{Lemma 2.8} Let $n\in\Bbb N$, $n\ge 10$ and let $T_n$ be a tree
with $n$ vertices and $\Delta(T_n)=n-4$. Suppose that for any
positive integer
 $m\ge n$ and connected
 graph $H\in \t{Ex}(m;T_n)$ we have $\Delta(H)\le n-4$.
Assume $p\in \Bbb N$, $p=k(n-1)+r\ge n-1$, where $k\in\Bbb N$ and
$r\in\{0,1,\ldots,n-2\}$. Then
$$\f{(n-2)p-r(n-1-r)}2\le \ex(p;T_n)\le \f{(n-2)p}2-\min\Big\{n-1+r,\f {r(n-1-r)}2\Big\}.$$
Hence, for $r\in\{0,1,2,n-5,n-4,n-3,n-2\}$ we have
$$\ex(p;T_n)=\f{(n-2)p-r(n-1-r)}2.$$
\endpro
Proof. Since $kK_{n-1}\cup K_r$ does not contain $T_n$ as a
subgraph, we see that
$$\ex(p;T_n)\ge e(kK_{n-1}\cup K_r)=\b{n-1}2+\b
r2=\f{(n-2)p-r(n-1-r)}2.$$ We claim that
$$\ex(n-1+r;T_n)\le
\f{(n-2)(n-1+r)}2-\min\Big\{n-1+r,\f{r(n-1-r)}2\Big\}.$$ As
$\ex(n-1;T_n)=e(K_{n-1})=\b{n-1}2$, we see that the claim holds for
$r=0$. Now suppose $r\ge 1$ and $G\in\t{Ex}(n-1+r;T_n)$. If $G$ is
connected, then $\Delta(G)\le n-4$ and so $e(G)\le \f
{(n-4)(n-1+r)}2=\f{(n-2)(n-1+r)}2-(n-1+r)$. Thus the claim is true.
\par Now suppose that $G$ is not connected and $G=G_1\cup
\cdots \cup G_s$, where $G_i$ is a component
 of $G$ with $|V(G_i)|=p_i$ and $p_1\le p_2\le \cdots\le p_s$. By
 Lemma 2.5, $p_1\le p_2=\cdots=p_{s-1}=n-1\le p_s$.
As $p_1+\cdots+p_s=n-1+r<2(n-1)$ we see that $s=2$, $p_1<n-1$ and
$p_2=n-1+r-p_1\ge n-1$. If $p_1>r$, then clearly $p_2<n-1$ and so
$e(G)=e(K_{p_1}\cup K_{n-1+r-p_1})=\b{p_1}2 +\b{n-1+r-p_1}2$. Using
Lemma 2.4(ii) we see that $$e(G)=\b{p_1}2
+\b{n-1+r-p_1}2<\b{n-1}2+\b r2=e(K_{n-1}\cup K_r).$$ This is a
contradiction. Hence $p_1\le r$.  If $p_1<r$, then $p_2=n-1+r-p_1\ge
n$ and so $\Delta(G)\le n-4$. Using Lemma 2.5 we see that $p_1\le
n-7$.  Hence
$$\align e(G)&=e(G_1)+e(G_2)\le \f{p_1(p_1-1)}2+\f{(n-4)(n-1+r-p_1)}2
\\&\q=\f{(n-4)(n-1+r)-p_1(n-3-p_1)}2
\\&\q<\f{(n-4)(n-1+r)}2=\f{(n-2)(n-1+r)}2-(n-1+r).\endalign$$ This shows
that the claim is also true for $p_1<r$. For $p_1=r$ we see that
$$\align e(G)&=e(K_{n-1}\cup K_r)=
\f{(n-1)(n-2)+r(r-1)}2\\&=\f{(n-2)(n-1+r)}2-\f {r(n-1-r)}2.
\endalign$$
 So the claim is also true. Hence the result is true for $p<2n-2$.
 \par Now assume $p\ge 2n-2$. By Lemma 2.7 and the above,
$$\align \ex(p;T_n)&=\f{(n-2)(p-(n-1+r))}2+\ex(n-1+r;T_n)
\\&\le \f{(n-2)(p-(n-1+r))}2+\f{(n-2)(n-1+r)}2
\\&\qq-\min\Big\{n-1+r,\f{r(n-1-r)}2\Big\}
\\&\q=\f{(n-2)p}2-\min\Big\{n-1+r,\f{r(n-1-r)}2\Big\}
.\endalign$$
To complete the proof, we note that $\f{r(n-1-r)}2\le
n-1-r$ for $r\in\{0,1,2,n-5,n-4,n-3,n-2\}$.

\pro{Lemma 2.9} Let $n\in\Bbb N$, $n\ge
 10$, $r\in\{0,1,\ldots,n-2\}$ and let $T_n$ be a tree with $n$ vertices
 and  $\Delta(T_n)=n-4$.
Suppose that for any positive integer
 $m\ge n$ and connected
 graph $H\in \t{Ex}(m;T_n)$ we have $\Delta(H)\le n-5$.
 Then
$$\ex(n-1+r;T_n)=\max\Big\{\big[\f {(n-5)(n-1+r)}2\big],\b
{n-1}2+\b r2\Big\}.$$
\endpro
Proof. Clearly $\ex(n-1;T_n)=e(K_{n-1})=\b{n-1}2$. Thus the result
is true for $r=0$. Now assume $r\ge 1$. By [3, Theorem 2.1],
$\ex(n-1+r;K_{1,n-4})=[\f
 {(n-5)(n-1+r)}2].$
 Since $\Delta(T_n)=n-4$ we see that
 $\ex(n-1+r;T_n)\ge \ex(n-1+r;K_{1,n-4})=[\f
 {(n-5)(n-1+r)}2].$ On the other hand,
 $\ex(n-1+r;T_n)\ge e(K_{n-1}\cup K_r)=\b{n-1}2+\b r2$.
 Thus,
$$\ex(n-1+r;T_n)\ge \max\Big\{\big[\f {(n-5)(n-1+r)}2\big],\b
{n-1}2+\b r2\Big\}.$$  Suppose $G\in  Ex(n-1+r;T_n).$ If $G$ is
connected, then $\Delta(G)\le n-5$ and so $e(G)\le \f
{(n-5)(n-1+r)}2$. Hence
$$\align \ex(n-1+r;T_n)&=e(G)\le \big[\f{(n-5)(n-1+r)}2\big]
\\&\q\le \max\Big\{\big[\f {(n-5)(n-1+r)}2\big],
\b {n-1}2+\b r2\Big\}.\endalign$$
This yields the result in this
case.

 \par Now suppose that $G$ is not connected and $G=G_1\cup
\cdots \cup G_s$, where $G_i$ is a component
 of $G$ with $|V(G_i)|=p_i$ and $p_1\le p_2\le \cdots\le p_s$.
By the argument in the proof of Lemma 2.8, we have $s=2$ and $p_1\le
r$.
 \par If $p_1<r$, then $p_2=n-1+r-p_1\ge n$. Using Lemma 2.5
we see that $p_1\le n-7$. By the assumption, $\Delta(G_2)\le n-5$
and so $e(G_2)\le [\f{(n-5)p_2}2]$. Hence
$$\align e(G)&=e(G_1)+e(G_2)\le \f{p_1(p_1-1)}2+\big[\f{(n-5)(n-1+r-p_1)}2
\big]\\&\q=\big[\f{(n-5)(n-1+r)-p_1(n-4-p_1)}2\big] \le
\big[\f{(n-5)(n-1+r)-3p_1}2\big]\\&\q<
\big[\f{(n-5)(n-1+r)}2\big].\endalign$$ This is a contradiction.
Thus, $p_1=r$ and so
$$e(G)=e(K_{n-1}\cup K_r)\le \max\Big\{\big[\f {(n-5)(n-1+r)}2\big],
\b {n-1}2+\b r2\Big\}.$$
\par  By the above, the lemma is proved.

\pro{Lemma 2.10} Let $n\in\Bbb N$, $n\ge 10$ and let $T_n$ be a tree
with $n$ vertices and $\Delta(T_n)=n-4$. Suppose that for any
positive integer
 $m\ge n$ and connected
 graph $H\in \t{Ex}(m;T_n)$ we have $\Delta(H)\le n-5$.
Assume $p=k(n-1)+r\ge n-1$, where $k\in\Bbb N$ and
$r\in\{0,1,\ldots,n-2\}$. Then
$$ex(p;T_n)=\f{(n-2)p-r(n-1-r)}2
+\max\Big\{0,\big[\f{r(n-4-r)-3(n-1)}2\big]\Big\}.$$
\endpro
 Proof. By Lemma 2.9,
$$\align &\ex(n-1+r;T_n)\\&=\max\Big\{\big[\f
 {(n-5)(n-1+r)}2\big],\b
{n-1}2+\b r2\Big\}
\\&=\f{(n-1)(n-2)+r(r-1)}2+
\max\Big\{0,\big[\f{r(n-4-r)-3(n-1)}2\big]\Big\}
\\&=\f{(n-2)(n-1+r)-r(n-1-r)}2+
\max\Big\{0,\big[\f{r(n-4-r)-3(n-1)}2\big]\Big\}.\endalign$$
 Thus the result is true for $p=n-1+r<2n-2$.
\par Now assume $p\ge 2n-2$.
From the above and Lemma 2.7 we see that
$$\align& \ex(p;T_n)\\&=\f{(n-2)(p-(n-1+r))}2+\max\Big\{\big[\f
 {(n-5)(n-1+r)}2\big],\b {n-1}2+\b r2\Big\}
 \\&=\max\Big\{\f{(n-2)p-r(n-1-r)}2,
 \big[\f{(n-2)p-3(n-1+r)}2\big]\Big\}
 \\ &=\f {(n-2)p-r(n-1-r)}2+\max\Big\{0,
 \big[\f{r(n-4-r)-3(n-1)}2\big]\Big\}.\endalign$$
This completes the proof.

\section*{3. Evaluation of $ex(p;T_n'')$}
\pro{Lemma 3.1} Let $p,n \in \Bbb N, p\geq n\geq 10$ and $G\in
Ex(p;T_n'')$. Suppose  that G is connected. Then $\Delta(G)=n-4$ or
$n-5$.\endpro
 Proof. By Lemma 2.1, $\Delta(G)\geq n-5$. Thus it is sufficient to prove
 that $\Delta(G)\le n-4$.
Suppose that $v_0\in V(G), d(v_0)=\Delta(G)=m$ and
$\Gamma(v_0)=\{v_1,\ldots,v_m\}$. If $p=m+1$, then
$V(G)=\{v_0,v_1,\ldots,v_m\}$ and $m=p-1\geq n-1$. Set
$G'=G[v_1,\ldots,v_m]$. If $d_{G'}(v_i)\ge 3$ for some
$i=1,2,\ldots,m$, as G does not contain $T_n''$ we see that
$e(G')=d_{G'}(v_i)\le m-1$. Otherwise, we have $d_{G'}(v_i)\le 2$
for every $i=1,2,\ldots,m$ and so $e(G')\le 2m/2=m$. Hence we always
have $$ e(G)=d(v_0)+e(G')\leq m+m =2p-2<\f{(n-5)p-1}2\le
\Big[\f{(n-5)p}2\Big].$$ This contradicts to (2.1). Thus $p>m+1$.

\par Suppose that
$u_1,\ldots,u_t$ are all vertices such that
$d(u_1,v_0)=\cdots=d(u_t,v_0)=2$ and
$\Gamma_2(v_0)=\{u_1,\ldots,u_t\}$. Then $t\geq 1$. Assume
$u_1v_1\in E(G)$ with no loss of generality.  Set
$V_1=\{v_0,v_1,\ldots,v_m\}$ and  $V_2=\{v_0,v_1,\ldots,v_m,u_1\}$.
\par Suppose $t=1$ and $m\ge n-2$.  If
$d_{G'}(v_1)\ge 3$, as $G$ does not contain $T_n''$ we see that
$\{v_2,v_3,\ldots,v_m\}$ is an independent set in $G'$. Hence
$e(G)-e(G-V_1)\le d(v_0)+d(u_1)+d(v_1)-2\le 3m-2.$  If
$d_{G'}(v_1)\le 2$, as $G$ does not contain $T_n$ we see that
$G[v_2,\ldots,v_m]$ does not contain $K_{1,2}$. Set
$G''=G[v_2,\ldots,v_m]$. Then $d_{G''}(v_i)\le 1$ for
$i=2,3,\ldots,m$ and so $e(G'')=\f
12\sum_{i=2}^md_{G''}(v_i)\le\f{m-1}2$. Therefore,
$$e(G)-e(G-V_1)=e(G[V_1])+d(u_1)
\le d(v_0)+2+e(G'')+d(u_1)\le m+2+\f{m-1}2+m<3m.$$ From the above
and Lemma 2.2 we see that $\Delta(G)\le n-3$ for $t=1$.
\par Suppose $t=1$ and $\Delta(G)=m\in\{n-3,n-4\}$. Then
$$\align e(G)-e(G-V_2)&\le d(u_1)+e(G[v_0,v_1,\ldots,v_m])
\\&\le m+e(K_{m+1})=\f{m^2+3m}2<\f{(m+1)(m+2)}2=e(K_{m+2}).
\endalign$$
Thus, $e(G)<e((G-V_2)\cup K_{m+2})$, which contradicts the
assumption $G\in Ex(p;T_n'')$.
\par By the above, for $t=1$ we have $\Delta(G)\le n-5$.
From now on we assume that $t\ge 2$.  Suppose that $m=\Delta(G)\ge
n-3$, $|\Gamma(v_1)\cap \Gamma_2(v_0)|\ge 2$ and
$u_1,u_2\in\Gamma(v_1)\cap \Gamma_2(v_0)$. Then
$\{v_2,v_3,\ldots,v_m\}$ is an independent set. If $t=2$ and
$v_1v_i\notin E(G)$ for $i=2,3,\ldots,m$, then $e(G)-e(G-V_1)\le
d(v_0)+d(u_1)+d(u_2)\le 3m.$  If $t=2$ and $v_1v_i\in E(G)$ for some
$i\in\{2,3,\ldots,m\}$, then $u_1v_j,u_2v_j\not\in E(G)$ for
$j\in\{2,3,\ldots,m\}-\{i\}$. Hence $e(G)-e(G-V_1) \le
d(v_1)+d(v_0)-1+2\le m+m+1.$ For $t\ge 3$ we must have
$u_3,\ldots,u_t\in\Gamma(v_1)$ and $u_iv_j\notin E(G)$ for
$i=1,2,\ldots,t$ and $j=2,3,\ldots,m$.
 Thus,
$e(G)-e(G-V_1) \le d(v_1)+d(v_0)-1\le m+m-1.$ From the above we
always have $e(G)-e(G-V_1)\le 3m$. This contradicts Lemma 2.2. Hence
$m=\Delta(G)\le n-4$ in the case $|\Gamma(v_1)\cap \Gamma_2(v_0)|\ge
2$. \par Now suppose that $u_1v_1,u_2v_2,\ldots,u_tv_t\in E(G)$ for
$t\ge 2$. We first assume $m=\Delta(G)\ge n-2$. If $t=2$ and
$v_1v_2\in E(G)$, then $d(v_3)=\cdots=d(v_m)=1$ and so
$$e(G)-e(G-V_1)\le d(v_1)+d(v_2)-1+d(v_3)+\cdots +d(v_m)\le
4+3+m-2=m+5<3m.$$ If $t=2$ and $v_1v_2\notin E(G)$, then clearly
$d(v_i)\le 3$ for $i=1,2,\ldots,m$. Hence $e(G)-e(G-V_1)\le
d(v_1)+d(v_2)+\cdots +d(v_m)\le 3m.$ For $t\ge 3$ we see that
$d(v_i)\le 2$ for $i=1,2,\ldots,m$. Thus, $e(G)-e(G-V_1)\le
d(v_1)+\cdots+d(v_m)\le 2m$. From the above we always have
$e(G)-e(G-V_1)\le 3m$, which contradicts Lemma 2.2. Therefore
$m=\Delta(G)\le n-3$.
\par Suppose $m=\Delta(G)=n-3$. If $t=2$, as G
does not contain any copies of $T_n''$ we see that $v_iv_j\notin
E(G)$ for $i\in \{1,2\}$ and $j\in \{3,4,\ldots,n-3\}$. Hence
$e(G[v_0,v_1,\ldots,v_{n-3}])\leq e(K_{n-2})-2(n-5).$ Thus,
$$\align e(G)-e(G-V_2)&\leq
d(u_1)+e(G[v_0,v_1,\ldots,v_{n-3}])+d(u_2)
\\&\leq \binom {n-2}2-2(n-5)+2(n-3)=\f {n^2-5n+14}2\\&<\f
{n^2-3n+2}2=e(K_{n-1}).
\endalign$$ This contradicts the fact $G\in
Ex(p;T_n'')$. If $t\geq 3$, then $\Gamma(v_i)=\{v_0,u_i\}$ for
$i=1,2,\ldots,t$. Hence
$$\align &e(G)-e(G-V_2)\\&\leq d(u_1)+1+d(v_2)+\cdots
+d(v_t)+e(G[v_0,v_{t+1},\ldots,v_{n-3}])+(t-1)(n-3-t)
\\&\leq n-3+1+2(t-1)+\f {(n-2-t)(n-3-t)}2+(t-1)(n-3-t)
\\&=\f {(n-1)(n-2)-(t^2-5t+2n-2)}2.
\endalign$$ For $t\geq 3$ we have $t^2-5t+2n-2\geq -6+2n-2>0$ and
so $e(G)-e(G-V_2)<\f {(n-1)(n-2)}2=e(K_{n-1}).$ That is,
$e(G)<e(K_{n-1}\cup (G-V_1))$, which contradicts the assumption.
Hence $\Delta(G)\le n-4$ as claimed.
 \pro{Lemma 3.2} Let $p,n \in
\Bbb N,\ p\geq n\geq 10$ and $G\in Ex(p;T_n'')$. Suppose  that G is
connected. Then $\Delta(G)=n-5$.\endpro
 Proof.  By Lemma 3.1, $m=\Delta(G)\le n-4$.
 Suppose that $v_0\in V(G), d(v_0)=\Delta(G)=n-4$,
$\Gamma(v_0)=\{v_1,\ldots,v_{n-4}\}$ and
$\Gamma_2(v_0)=\{u_1,\ldots,u_t\}$. By the proof of Lemma 3.1, we
have $\Delta(G)=n-5$ for $t=1$. From now on we assume $t\ge 2$. If
$t=2$ and $u_1v_1,u_2v_1\in E(G)$, setting
$V_1=\{v_0,v_1,\ldots,v_{n-4}\}$ and
$V_2=\{v_0,v_1,\ldots,v_{n-4},u_1,u_2\}$ we see that
$$\align e(G)-e(G-V_2)&\le d(u_1)+d(u_2)+e(G[V_1])\\&\le n-4+n-4+\f{(n-3)(n-4)}2
=\f{n^2-3n-4}2\\&<\f{(n-1)(n-2)}2=e(K_{n-1})\endalign$$ and so
$e(G)<e((G-V_2)\cup K_{n-1})$. This contradicts the assumption $G\in
Ex(p;T_n'')$. If $t\ge 3$ and $v_1u_i\in E(G)$ for $i=1,2,\ldots,t$,
then $u_iv_j\not\in E(G)$ for $i=1,2,\ldots,t$ and
$j=2,3,\ldots,n-4$. Thus,
$$\align e(G)-e(G-V_2)&\le d(u_1)+d(u_2)+d(v_1)
+e(G[v_0,v_2,v_3,\ldots,v_{n-4}])\\&\le n-4+n-4+n-4+\f{(n-4)(n-5)}2
=\f{n^2-3n-4}2\\&<\f{(n-1)(n-2)}2=e(K_{n-1})\endalign$$ and so
$e(G)<e((G-V_2)\cup K_{n-1})$, which contradicts the assumption
$G\in Ex(p;T_n'')$.

\par Now suppose that $u_1v_1,\ldots,u_tv_t\in E(G)$.
 Then $2\leq t\leq n-4$. As $d(v_1)\le n-4$ we see that
 $v_1v_i\not\in E(G)$ for some $i\in\{2,3,\ldots,n-4\}$.
Thus, for $t=2$ we have
$$\align &e(G)-e(G-V_2)\\&
\le d(u_1)+d(u_2)+e(G[v_0,v_1,\ldots,v_{n-4}])
\\&\le n-4+n-4+\b{n-3}2-1=\f{n^2-3n-6}2<\b{n-1}2=e(K_{n-1}).
\endalign$$
This yields $e(G)<e((G-V_2)\cup K_{n-1})$, which is impossible.
Hence $t\ge 3$.
\par Now suppose $u_1v_1,\ldots,u_tv_t\in E(G)$, $t\ge 3$ and
 $V_1=\{v_0,v_1,\ldots,v_{n-4},u_1\}$.
 We first
claim that $d(v_i)\leq n-5$ for $i=1,2,\ldots,t$. Suppose
$d(v_1)=n-4$. Then $v_1v_i\not\in E(G)$ for some
$i\in\{2,3,\ldots,n-4\}$ and so $|\Gamma(u_1)\cap(G-V_1)|\le 1$.
Otherwise, for $w_1,w_2\in \Gamma(u_1)\cap(G-V_1)$,
$G[v_1,v_2,\ldots,v_{n-4},v_0,u_1,w_1,w_2]$ contains a copy of
$T_n''$. For $i\in\{1,2,\ldots,n-4\}$ there is at most one vertex in
$\{u_1,\ldots,u_t\}$ adjacent to $v_i$. Hence
$$\align &e(G)-e(G-V_1)\\&\le e(G[v_0,v_1,\ldots,v_{n-4}])+n-4+1
\\&\le \b{n-3}2-1+n-4+1=\f{n^2-5n+4}2<\b{n-2}2=e(K_{n-2})
\endalign$$
and so $e(G)<e((G-V_1)\cup K_{n-2})$. This is impossible. Hence the
claim is true. Set $G'=G[v_0,v_1,\ldots,v_{n-4}]$. Suppose that
there are exactly $s$ vertices $v_{i_1},\ldots,v_{i_s}$ in
$\{v_2,\ldots,v_{n-4}\}$ adjacent to some vertex in
$\{u_2,\ldots,u_t\}$. Then $d(v_{i_j})\le n-5$ for $j=1,2,\ldots,s$
by the above argument. Hence
$$e(G')=\f 12\sum_{i=0}^{n-4}d_{G'}(v_i)\leq \f {(s+1)(n-6)+(n-4-s)(n-4)}2=\f
{n^2-7n+10-2s}2$$ and therefore
$$\align e(G)-e(G-V_1)
&\leq e(G')+d(u_1)+s\leq \f {n^2-7n+10}2-s+n-4+s\\&=\f
{n^2-5n+2}2<\f {(n-2)(n-3)}2=e(K_{n-2}).\endalign$$ This contradicts
the fact $G\in Ex(p;T_n'')$. Thus $\Delta(G)=n-5$ as claimed.

\pro{Theorem 3.1} Let $p,n \in \Bbb N, p\geq n\geq 10$,
$p=k(n-1)+r$, $k\in\Bbb N$ and $r\in\{0,1,\ldots,n-2\}$. Then
$$ex(p;T_n'')=\f{(n-2)p-r(n-1-r)}2+\max\Big\{0,\big[\f{r(n-4-r)-3(n-1)}2\big]\Big\}.$$\endpro
Proof. This is immediate from Lemmas 3.2 and 2.10.

\section*{4. Evaluation of $ex(p;T_n^3)$}

 \pro{Lemma 4.1} Let $p,n \in
\Bbb N, p\geq n\geq 10$ and $G\in Ex(p;T_n^3)$. Suppose that G is
connected. Then $\Delta(G)=n-5$ or $n-4$.\endpro
 Proof. Suppose that
$v_0\in V(G), d(v_0)=\Delta(G)=m$ and
$\Gamma(v_0)=\{v_1,\cdots,v_m\}$. If $p=m+1$, then
$V(G)=\{v_0,v_1,\cdots,v_m\}$ and $m=p-1\geq n-4$. Since G does not
contain $T_n^3$, we see that $G{[v_1,v_2,\cdots,v_m]}$ does not
contain $K_{1,3}$ and hence $\Delta(G{[v_1,v_2,\cdots,v_m]})\le 2$.
Thus, $e(G{[v_1,v_2,\cdots,v_m]})\leq m$. Therefore $$\aligned
e(G)=d(v_0)+e(G{[v_1,v_2,\cdots,v_m]})\leq m+m =2p-2< [\f{(n-5)p}2].
\endaligned$$
This contradicts to (2.1). Thus $p>m+1$.

\par Suppose that $\Gamma_2(v_0)=\{u_1,\cdots,u_t\}$. Then $t\geq 1$. We may suppose
that $ v_1,\cdots,v_{s_1} $ are all vertices adjacent to exactly two
vertices in the set $\{u_1,\cdots,u_t\}$ and
$v_{s_1+1},\cdots,v_{s_2}$ are all vertices adjacent to exactly one
vertex in the set $\{u_1,\cdots,u_t\}$. Let $
V_1=\{v_0,v_1,\cdots,v_m\}$, $V_1'=V(G)-V_1$ and let $e(V_1V_1')$ be
the number of edges with one endpoint in $V_1$ and another endpoint
in $V_1'$. As $G$ does not contain $T_n^3$ we see that
$e(V_1V_1')=2s_1+s_2-s_1=s_1+s_2$.
\par
If $m\geq n-1$, as G does not contain $T_n^3$ as a subgraph, we see
that $d(v_i)\le 3$ for $i=1,2,\ldots,m$ and so $e(G)\le
3m+e(G-V_1)$. This contradicts Lemma 2.2. Hence $m\le n-2$.
\par
Suppose $m=n-2$. Since $G$ does not contain $T_n^3$ we see that
$d(v_i)\le 3$ for $i=1,2,\ldots,s_2$. Thus,
$$e(G)-e(G-V_1)\le 3s_2+\b{n-1-s_2}2=s_2(s_2-(2n-6))+\b{n-1}2.$$
As $s_2\le n-2<2n-6$ we have $e(G)<\b {n-1}2+e(G-V_1)=e(K_{n-1}\cup
(G-V_1)).$ This contradicts the assumption $G\in Ex(p; T_n^3)$.
Hence $m\le n-3$.
\par
Suppose $m=n-3$. For $i\in\{1,2,\ldots,s_1\}$ and
$j\in\{s_1+1,\ldots,n-3\}$ we have $v_iv_j\not\in E(G)$. Thus,
$$\aligned e(G)-e(G-V_1)&=n-3+2s_1+(s_2-s_1)+e(G[v_{s_1+1},\ldots,v_{n-3}])
\\&\le n-3+2s_1+s_2-s_1+\b{n-3-s_1}2
\\&=\b{n-2}2+s_2-\f{s_1(2n-9-s_1)}2.\endaligned$$
If $s_1\ge 2$, then $\f{s_1(2n-9-s_1)}2\ge 2n-11>n-3\ge s_2$ and
hence $e(G)<e(G-V_1)+\b{n-2}2=e((G-V_1)\cup K_{n-2})$. This
contradicts the fact $G\in Ex(p; T_n^3)$. Hence $s_1=0$ or $1$. We
claim that $d(v_i)\ge n-4$ for some $i\in\{1,2,\ldots,s_2\}$. Assume
that $d(v_i)\le n-5$ for $i=1,2,\ldots,s_2$. If $s_1=0$, then
$e(G[V_1])\le \f{s_2(n-6)+(n-2-s_2)(n-3)}2=\b{n-2}2-\f 32s_2$ and so
$$e(G)\le \b{n-2}2-\f 32s_2+s_2+e(G-V_1)<e(K_{n-2})+e(G-V_1)
=e(K_{n-2}\cup (G-V_1)).$$ This contradicts the fact $G\in Ex(p;
T_n^3)$. If $s_1=1$, we may assume that there are two vertices in
$G-V_1$ adjacent to $v_1$. Then $d(v_1)=3$ and $v_1v_i\not\in E(G)$
for $i=2,3,\ldots,n-3$. Thus,
$$e(G[V_1])\le
\f{1+(s_2-1)(n-6)+(n-2-s_2)(n-3)}2=\f{(n-2)(n-3)-3s_2-(n-7)}2$$ and
so
$$\aligned e(G)&\le e(G-V_1)+s_2+1+\f{(n-2)(n-3)-3s_2-(n-7)}2
\\&=e(G-V_1)+\b{n-2}2-\f{s_2+n-9}2<e((G-V_1)\cup K_{n-2}).\endaligned$$ This
contradicts the fact $G\in Ex(p; T_n^3)$. Hence the claim is true.
Now suppose $d(v_1)\ge n-4$ and $u_1v_1\in E(G)$, where $u_1\in
G-V_1$. Set $V_2=\{u_1,v_0,v_1,\ldots,v_{n-3}\}$ and
$V_2'=V(G)-V_2$. Then there are at most two vertices in $V_2'$
adjacent to $u_1$. Suppose that there are exactly $r$ vertices in
$\{v_1,\ldots,v_{s_2}\}$ adjacent to $u_1$. As $s_1=0$ or $1$ we
have $e(V_2V_2')\le s_2-r+1+2$. As $\Delta(G)\le n-3$ we see that
$d_{G[V_2]}(v_i)\le n-4$ for $i=1,2,\ldots,s_2$. Thus,
$$e(G[V_2])=r+e(G[V_1])\le r+\f{s_2(n-4)+(n-2-s_2)(n-3)}2
=\b{n-2}2-\f{s_2}2+r.$$ Therefore,
$$\aligned e(G)&=e(G[V_2])+e(V_2V_2')+e(G-V_2)
\\&\le \b{n-2}2-\f {s_2}2+r+s_2-r+3+e(G-V_2)
=\b{n-2}2+\f{s_2+6}2+e(G-V_2)
\\&\le \f{(n-2)(n-3)+n-3+6}2+e(G-V_2)
\\&<\b{n-1}2+e(G-V_2)=e(K_{n-1}\cup(G-V_2)),
\endaligned$$
which contradicts the assumption $G\in Ex(p; T_n^3)$. Hence
$\Delta(G)\le n-4$.

\pro{Theorem 4.1} Let $p,n \in \Bbb N,\ p\ge n\ge 10$, $p=k(n-1)+r$,
$k\in\Bbb N$ and $r\in\{0,1,2,n-5,n-4,n-3,n-2\}$. Then
$$\ex(p;T_n^3)=\f{(n-2)p-r(n-1-r)}2.$$
\endpro
Proof. This is immediate from Lemmas 4.1 and 2.10.

\pro{Lemma 4.2} Let $n \in \Bbb N,\  n\geq 10$, $p=n-1+r$,
$r\in\{1,2,\ldots,n-6\}$ and $G\in Ex(p;T_n^3)$. Suppose that G is
connected. If $r(n-8-r)>5+((-1)^n-(-1)^{(n-1)(r-1)})/2$, then
$\Delta(G)=n-5$.
\endpro
Proof. By Lemma 4.1, $\Delta(G)=n-4$ or $n-5$.
 Suppose
$\Delta(G)=n-4$, $v_0\in V(G)$, $d(v_0)=n-4$,
$\Gamma(v_0)=\{v_1,\cdots,v_{n-4}\}$ and
$\Gamma_2(v_0)=\{u_1,\cdots,u_t\}$. Then $1\le t\le n-4$. We may
suppose that $ v_1,\cdots,v_{s_1} $ are all vertices adjacent to
exactly two vertices in $\Gamma_2(v_0)$ and
$v_{s_1+1},\cdots,v_{s_2}$ are all vertices adjacent to exactly one
vertex in $\Gamma_2(v_0)$.
\par
We first claim that $d(v_i)\leq n-5$ for $i=1,2,\ldots,s_1$. Suppose
$d(v_1)=n-4$ and $u_1$,$u_2$ are two vertices in
$G-\{v_0,v_1,\ldots,v_{n-4}\}$ adjacent to $v_1$. Set
$G'=G[v_0,v_1,\ldots,v_{n-4}]$. Then $$\aligned e(G')&=\f
12(d_{G'}(v_0)+\sum_{i=1}^{s_1}d_{G'}(v_i)+\sum _{i=s_1+1}^{s_2}
d_{G'}(v_i)+\sum_{i=s_2+1}^{n-4} d_{G'}(v_i))
\\&\leq \f 12(n-4+(n-6){s_1}+(n-5)(s_2-s_1)+(n-4-s_2)(n-4))
\\&=\f {n^2-7n+12-s_1-s_2}2
\endaligned$$
Let $V_1=\{v_0,v_1,\ldots,v_{n-4},u_1,u_2\}$ and $V_1'=V(G)-V_1$. As
$d(v_1)=n-4$, for $i=1,2$ there are at most two vertices in $G-V_1$
adjacent to $u_i$. Thus, $$\align
e(G)-e(G-V_1)&=e(G[V_1])+e(V_1V_1')\leq e(G')+2s_1+s_2-s_1+2+2
\\&\leq \f {n^2-7n+12-s_1-s_2}2+s_1+s_2+4=\f {n^2-7n+20+s_1+s_2}2
\\&\leq \f {n^2-7n+20+n-4+n-4}2=\f {n^2-5n+12}2
\\&<\f {(n-1)(n-2)}2=e(K_{n-1}).
\endalign$$
Hence $ex(p;T_n^3)=e(G)<e(G-V_1)+e(K_{n-1})=e((G-V_1)\cup K_{n-1})$.
As $(G-V_1)\cup K_{n-1}$ does not contain $T_n^3$, we get a
contradiction. Hence $d(v_i)\leq n-5$ for $i=1,2,\ldots,s_1$.
\par Now we show that $d(v_i)\leq n-5$ for $s_1<i\leq s_2$. Suppose
$d(v_i)=n-4$ for $i\in \{s_1+1,\ldots,s_2\}$ and $u_iv_i\in E(G)$,
where $u_i\in V(G)-\{v_0,v_1,\ldots,v_{n-4}\}$. For
$G'=G[v_0,v_1,\ldots,v_{n-4}]$, from the above and the fact
$\Delta(G)\leq n-4$ we see that $d_{G'}(v_i)\leq n-7$ for $1\leq
i\leq s_1$, $d_{G'}(v_i)\leq n-5$ for $s_1< i\leq s_2$ and
$d_{G'}(v_i)\leq n-4$ for $s_2< i\leq n-4$. Thus
$$\align e(G')&=\f 12\sum_{i=0}^{n-4}d_{G'}(v_i)\\&\leq \f
12(n-4+(n-7)s_1+(n-5)(s_2-s_1)+(n-4-s_2)(n-4))\\&\q=\f{n^2-7n+12-2s_1-s_2}2.
\endalign$$
Set $V_1=\{v_0,v_1,\ldots,v_{n-4},u_1\}$ and $V_1'=V(G)-V_1$. As
$d(v_1)=n-4$, there are at most two vertices in $G-V_1$ adjacent to
$u_1$. Note that $s_2\leq n-4$. We deduce that $$\align
&e(G)-e(G-V_1)=e(G[V_1])+e(V_1V_1')\\&\q\le \f
{n^2-7n+12-2s_1-s_2}2+2s_1+s_2-s_1+2=\f {n^2-7n+16+s_2}2\\&\q\le \f
{n^2-7n+16+n-4}2=\f {n^2-6n+12}2 <\f
{(n-2)(n-3)}2=e(K_{n-2}).\endalign$$ Hence $e(G)<e((G-V_1)\cup
K_{n-2})$. This contradicts the assumption $G\in Ex(p;T_n^3)$.
\par By the above,
$$d(v_i)\le n-5\qtq{for} i=1,2,\ldots,s_2.\tag 4.1$$
 For $V_1=\{v_0,v_1,\ldots,v_{n-4}\}$ we have
$|V(G-V_1)|=p-(n-3)=r+2<n$ and so $e(G-V_1)\leq \binom {r+2}2$.  As
$\Delta(G)\leq n-4$ and $d_G(v_i)\leq n-5$ for $i=1,2,\ldots,s_2$,
we see that $d_{G[V_1]}(v_i)\leq n-7$ for $1\leq i\leq s_1$,
$d_{G[V_1]}(v_i)\leq n-6$ for $s_1<i\leq s_2$ and
$d_{G[V_1]}(v_i)\leq n-4$ for $s_2<i\leq n-4$. Thus,
$$\aligned e(G[V_1])&=\f 12\sum_{i=0}^{n-4}d_{G[V_1]}(v_i)
 \\&\le \f 12(n-4+(n-7)s_1+(n-6)(s_2-s_1)+(n-4)(n-4-s_2))
\\&\q=\f{n^2-7n+12-s_1-2s_2}2.\endaligned\tag 4.2$$
Set $V_1'=V(G)-V_1$. Then
 $$\align
e(G)&=e(G[V_1])+e(V_1V'_1)+e(G-V_1)
\\&\leq \f
{n^2-7n+12-s_1-2s_2}2+2s_1+s_2-s_1+\binom {r+2}2\\&\q=\f
{n^2-7n+12+s_1+(r+1)(r+2)}2
\\&\leq \f {n^2-7n+12+n-4+r^2+3r+2}2=\f {n^2-6n+10+r^2+3r}2
\endalign$$
and so
$$e(G)\le \Big[\f {n^2-6n+10+r^2+3r}2\Big]=\f {n^2-6n+10+r^2+3r
-(1-(-1)^n)/2}2.$$ Suppose $G_0\in Ex(n-1+r;K_{1,n-4})$. Then
$$e(G_0)=\Big[\f{(n-1+r)(n-5)}2\Big]=\f {(n-1+r)(n-5)-
(1-(-1)^{(n-1)(r-1)})/2}2.$$ As $G_0$ does not contain $T_n^3$ and
$G\in Ex(n-1+r;T_n^3)$, we get $$\align &\f
{n^2-6n+10+r^2+3r-(1-(-1)^n)/2}2\\&\geq e(G)\geq e(G_0)\geq \f
{(n-1+r)(n-5)- (1-(-1)^{(n-1)(r-1)})/2}2
 \endalign$$ and so
$r(n-8-r)\leq 5+((-1)^n-(-1)^{(n-1)(r-1)})/2$. Hence, if
$r(n-8-r)>5+((-1)^n-(-1)^{(n-1)(r-1)})/2$, we must have $\Delta
(G)<n-4$ and so $\Delta (G)=n-5$ as claimed.

\pro{Lemma 4.3} Let $n,r\in\Bbb N$, $n\ge 15$ and $3\le r\le n-9$.
Then
$$\ex(n-1+r;T_n^3)=\max\Big\{\Big[\f{(n-5)(n-1+r)}2\Big],\b{n-1}2+\b
r2\Big\}.$$
\endpro
Proof. By the proof of Lemma 2.9, we have
$$\ex(n-1+r;T_n^3)\ge \max\Big\{\big[\f{(n-5)(n-1+r)}2\big],\b{n-1}2+\b
r2\Big\}.$$ For $3\le r \le n-10$ we see that $r(n-8-r)\ge
2(n-10)>7>5+((-1)^n-(-1)^{(n-1)(r-1)})/2$. For $r=n-9$ we also have
$r(n-8-r)=n-9>5+((-1)^n-(-1)^{(n-1)(r-1)})/2$. Let
$G\in\t{Ex}(n-1+r;T_n^3)$. If $G$ is connected, by Lemma 4.2 we have
$\Delta(G)\le n-5$ and hence $e(G)\le [\f{(n-5)(n-1+r)}2]$.
 \par Now suppose that $G$ is not connected and $G=G_1\cup
\cdots \cup G_s$, where $G_i$ is a component
 of $G$ with $|V(G_i)|=p_i$ and $p_1\le p_2\le \cdots\le p_s$.
By Lemma 4.1 and the argument in the proof of Lemma 2.8, we have
$s=2$ and $p_1\le r$. Set $r'=r-p_1$. Then $0\le r'=r-p_1\le
n-9-p_1\le n-10$. For $r'\ge 2$ we have $r'(n-8-r')\ge 2(n-10)>7>
5+((-1)^n-(-1)^{(n-1)(r'-1)})/2$. For $r'=1$ we have $r'(n-8-r')
=n-9>5+((-1)^n-(-1)^{(n-1)(r'-1)})/2$. Since
$|V(G_2)|=p_2=n-1+r-p_1=n-1+r'$, using Lemma 4.2 we see that for
$r'\ge 1$ we have $\Delta(G_2)\le n-5$ and so $e(G_2)\le
[\f{(n-5)(n-1+r-p_1)}2]$. From Lemmas 4.1 and 2.5 we see that
$p_1\le n-7$.  Hence, for $p_1<r$,
$$\align e(G)&=e(G_1)+e(G_2)\le \f{p_1(p_1-1)}2+\big[\f{(n-5)(n-1+r-p_1)}2
\big]\\&\q=\big[\f{(n-5)(n-1+r)-p_1(n-4-p_1)}2\big] \le
\big[\f{(n-5)(n-1+r)-3p_1}2\big]\\&\q<
\big[\f{(n-5)(n-1+r)}2\big].\endalign$$ This contradicts the fact
that $e(G)=\ex(n-1+r;T_n^3)\ge [\f{(n-5)(n-1+r)}2]$. Thus, $p_1=r$
and so $e(G)=e(K_{n-1}\cup K_r)=\b{n-1}2+\b r2$.
\par By the above, we always have
$$\ex(n-1+r;T_n^3)=e(G)\le \max\Big\{\Big[\f {(n-5)(n-1+r)}2\Big],
 \b {n-1}2+\b r2\Big\}.$$
 Thus the result is true.
\pro{Theorem 4.2} Let $p,n\in\Bbb N$, $p\ge n\ge 15$, $p=k(n-1)+r$,
 $k\in\Bbb N$ and $r\in\{3,4,\ldots,n-9\}$. Then
$$ex(p;T_n^3)=\f{(n-2)p-r(n-1-r)}2
+\max\Big\{0,\big[\f{r(n-4-r)-3(n-1)}2\big]\Big\}.$$
\endpro
 Proof. By Lemma 4.3,
$$\align &\ex(n-1+r;T_n^3)\\&=\max\Big\{\big[\f
 {(n-5)(n-1+r)}2\big],\b
{n-1}2+\b r2\Big\}
\\&=\f{(n-2)(n-1+r)-r(n-1-r)}2+
\max\Big\{0,\big[\f{r(n-4-r)-3(n-1)}2\big]\Big\}.\endalign$$
 Thus the result is true for $p=n-1+r<2n-2$.
\par Now assume $p\ge 2n-2$.
From the above and Lemmas 4.1 and 2.7 we see that
$$\align& \ex(p;T_n^3)\\&=\f{(n-2)(p-(n-1+r))}2+\ex(n-1+r;T_n^3)
\\&=\f{(n-2)(p-(n-1+r))}2+\max\Big\{\big[\f
 {(n-5)(n-1+r)}2\big],\b {n-1}2+\b r2\Big\}
  \\ &=\f {(n-2)p-r(n-1-r)}2+\max\Big\{0,
 \big[\f{r(n-4-r)-3(n-1)}2\big]\Big\}.\endalign$$
This completes the proof.

\pro{Lemma 4.4} Let $m,n \in \Bbb N$ with $m\le n-4$ and $n\ge 10$.
Suppose that $G\in Ex(2n-6-m;T_n^3)$ and $G$ is connected. Assume
that $v_0\in V(G)$ and $d(v_0)=\Delta(G)=n-4$. Then for any $v\in
V(G)- \{v_0\}\cup \Gamma(v_0)$ we have $d(v)\le n-5$.\endpro
 Proof. Assume that
$\Gamma(v_0)=\{v_1,\ldots,v_{n-4}\}$ and
$\Gamma_2(v_0)=\{u_1,\cdots,u_t\}$. Clearly $t\le n-3-m\le n-4$ and
$d(v)\le n-4-m\le n-5$ for $v\in
V(G)-\{v_0,v_1,\ldots,v_{n-4},u_1,\ldots,u_t\}$. Thus, we only need
to prove that $d(u_i)\le n-5$ for $i=1,2,\ldots,t$. We may suppose
that $v_1,\cdots,v_{s_1} $ are all vertices adjacent to exactly two
vertices in the set $\{u_1,\cdots,u_t\}$ and
$v_{s_1+1},\cdots,v_{s_2}$ are all vertices adjacent to exactly one
vertex in the set $\{u_1,\cdots,u_t\}$. Let
$V_1=\{v_0,v_1,\ldots,v_{n-4}\}$ and $V_1'=V(G)-V_1$.   Since
$n-5-m\le n-6$, by (4.2) we have $e(G[V_1])\le \f
{n^2-7n+12-s_1-2s_2}2$.  Suppose $d(u_i)=n-4$ for some
$i\in\{1,\ldots,t\}$ and
$\Gamma(u_i)\cap\{v_1,\ldots,v_{n-4}\}=\{v_{i_1},v_{i_2},\ldots,v_{i_k}\}$.
Then $m\le k\le n-4$. If $k\ge \f {n-8}2$, then $d_{G-V_1}(u_i)\le
n-4-\f {n-8}2=\f n2$ and so
$$\align e(G-V_1)&\le d_{G-V_1}(u_i)+\b{n-4-m}2\\&\le
\f n2+\f{(n-4-m)(n-5-m)}2= \f{n^2-(2m+8)n+(m+4)(m+5)}2. \endalign$$
Therefore,
$$\align &e(G)=e(G[V_1])+e(V_1V'_1)+e(G-V_1)
\\&\leq \f
{n^2-7n+12-s_1-2s_2}2+2s_1+s_2-s_1+\f{n^2-(2m+8)n+(m+4)(m+5)}2
\\&=\f{2n^2-(2m+15)n+m^2+9m+32+s_1}2
\\&\le \f{2n^2-(2m+15)n+m^2+9m+32+n-4}2=n^2-(m+7)n+\f{m(m+9)}2+14
\\&<n^2-(m+7)n+\f{m(m+11)}2+16=e(K_{n-1}\cup K_{n-5-m}) .\endalign$$ This is a
contradiction. Hence $k<\f {n-8}2$. As $d(u_i)=n-4$ and $G$ does not
contain $T_n^3$ as a subgraph, for $j\in\{i_1,\ldots,i_k\}$ we see
that $|\Gamma(v_j)\cap
(\{v_1,\ldots,v_{n-4}\}-\{v_{i_1},\ldots,v_{i_k}\})|\le 1$. Hence
$d_{G[V_1]}(v_j)\le k-1+1+1=k+1$. As $d(v_0)=n-4$ and $v_0v_j\in
E(G)$ we have $|\Gamma(v_j)\cap\{u_1,\ldots,u_t\}|\le 2$. Thus,
  $d(v_j)\le k+3$ and so
$d(v_{i_1})+\cdots +d(v_{i_k})\le k(k+3)$. For $m\le k< \f {n-8}2$
we see that
$$k(n-7-k)-m(n-7-m)=(k-m)(n-7-k-m)\ge (k-m)(n-7-\f{n-7}2-m)\ge 0.$$
 Thus, $k(n-7-k)\ge m(n-7-m)$ and so
$$\align e(G)=\f 12\sum_{v\in V(G)}d(v)&\le
\f{(2n-6-m-k)(n-4)+k(k+3)}2
\\&=\f{(2n-6-m)(n-4)-k(n-7-k)}2\\&\le \f{(2n-6-m)(n-4)-m(n-7-m)}2
\\&=n^2-(m+7)n+\f{m(m+11)}2+12\\&<n^2-(m+7)n+\f{m(m+11)}2+16
=e(K_{n-1}\cup K_{n-5-m}) .\endalign$$
 This contradicts the assumption $G\in Ex(2n-6-m;T_n^3)$. Hence $d(u_i)\le
n-5$ for $i=1,2,\ldots,t$ as claimed. The proof is now complete.

\pro{Lemma 4.5} Let $n \in \Bbb N$ with $n\ge 10$. Then
$ex(2n-7;T_n^3)=n^2-8n+22$.

\endpro
Proof. Let $G\in Ex(2n-7;T_n^3)$. As $K_{n-1}\cup K_{n-6}$ does not
contain $T_n^3$ as a subgraph, we have $e(G)\ge e(K_{n-1}\cup
K_{n-6})$. We first assume that $G$ is connected. By Lemma 4.1,
$\Delta(G)=n-5$ or $n-4$. If $\Delta (G)=n-5$, then
$$e(G)\le \f {(n-5)(2n-7)}2=n^2-\f {17}2n+\f
{35}2<n^2-8n+22=e(K_{n-1}\cup K_{n-6}),$$ which contradicts the fact
$e(G)\ge e(K_{n-1}\cup K_{n-6})$. Hence $\Delta (G)=n-4$. Suppose
that $v_0\in V(G)$, $d(v_0)=n-4$ and
$\Gamma(v_0)=\{v_1,\ldots,v_{n-4}\}$. By Lemma 4.4, $d(v)\le n-5$
for $v\in V(G)-\{v_0,v_1,\ldots,v_{n-4}\}$. Therefore,
$$\align
e(G)&=\f 12\sum_{v\in V(G)}d(v)\le\f{ (n-3)(n-4)+(n-4)(n-5)}2 \\&\q=
n^2-8n+16<n^2-8n+22=e(K_{n-1}\cup K_{n-6}).\endalign$$ This is also
a contradiction. So $G$ is not connected.
\par
Suppose that $G$ is not connected and $G=G_1\cup \cdots \cup G_s$,
where $G_i$ is a component
 of $G$ with $|V(G_i)|=p_i$ and $p_1\le p_2\le \cdots\le p_s$.
By Lemmas 4.1, 2.5 and the argument in the proof of Lemma 2.8, we
have $s=2$ and $p_1\le n-6$. If $p_1\le n-7$, then $p_2=2n-7-p_1\ge
n$. By Lemmas 4.1 and 2.5, we  have $2n-7=p_1+p_2\le 2n-8$. This is
impossible. Hence $p_1=n-6$, $p_2=n-1$ and so $e(G)=e(K_{n-1}\cup
K_{n-6})=n^2-8n+22$. This proves the lemma.
 \pro{Theorem 4.3} Let $p,n\in\Bbb N$, $p\ge n\ge 10$ and
$p=k(n-1)+n-6$ with $k\in\Bbb N$. Then
$$ex(p;T_n^3)=\f{(n-2)p-5(n-6)}2.$$
\endpro
Proof.  By Lemmas 4.1, 2.7 and 4.5,
$$\align \ex(p;T_n^3)&=\f{(n-2)(p-(2n-7))}2+\ex(2n-7;T_n)
\\&=\f{(n-2)(p-(2n-7))}2+n^2-8n+22
=\f{(n-2)p-5(n-6)}2.\endalign$$

\pro{Lemma 4.6} Let $n \in \Bbb N$ with $n\ge 15$. Then
$$ex(2n-9;T_n^3)=n^2-10n+24+\max\Big\{\big[\f n2\big],13\Big\}.$$
\endpro
Proof. Let $G\in Ex(2n-9;T_n^3)$. Suppose that $G$ is not connected
and $G=G_1\cup \cdots \cup G_s$, where $G_i$ is a component
 of $G$ with $|V(G_i)|=p_i$ and $p_1\le p_2\le \cdots\le p_s$.
By Lemmas 4.1, 2.5 and the argument in the proof of Lemma 2.8, we
have $s=2$ and $p_1\le n-8$. If $p_1\le n-9$, then $p_2=2n-9-p_1\ge
n$. By Lemmas 4.1 and 2.5, we  have $p_1(n-3-p_1)\le
p_1+p_2+1=2n-8$. For $3\le p_1\le n-9$ we have $p_1(n-3-p_1)\ge
3(n-6)>2n-8$. Thus, $p_1=1$ or $2$. For $p_1=1$ we have
$p_2=2n-10=n-1+n-9$. As
$(n-9)(n-8-(n-9))=n-9>5+((-1)^n-(-1)^{(n-1)(n-10)})/2$, using Lemma
4.2 we see that $\Delta(G_2)\le n-5$ and hence
$e(G)=e(G_2)=\ex(2n-10;K_{1,n-4})=\f{(2n-10)(n-5)}2=n^2-10n+25$. For
$p_1=2$ we have $p_2=2n-11=n-1+n-10$. As
$(n-10)(n-8-(n-10))=2(n-10)>7>5+((-1)^n-(-1)^{(n-1)(n-11)})/2$,
using Lemma 4.2 we see that $\Delta(G_2)\le n-5$ and hence
$e(G_2)=\ex(2n-11;K_{1,n-4})=[\f{(2n-11)(n-5)}2]=[\f{2n^2-21n+55}2]$.
Thus, $e(G)=e(G_1)+2(G_2)=1+[\f{2n^2-21n+55}2]=[\f{2n^2-21n+57}2]$.
For $p_1=n-8$ we have $p_2=n-1$ and so $e(G)=e(K_{n-1}\cup K_{n-8})
=\b{n-1}2+\b{n-8}2=n^2-10n+37$. Therefore, when $G$ is not
connected, we have
$$e(G)=\max\Big\{n^2-10n+25,\big[\f{2n^2-21n+57}2\big],n^2-10n+37\Big\}
=n^2-10n+37.$$

\par Assume that $G$ is connected.
By Lemma 4.1, $\Delta(G)\le n-4$. If $\Delta(G)\le n-5$, then
clearly
$e(G)=\ex(2n-9;K_{1,n-4})=[\f{(2n-9)(n-5)}2]=[\f{2n^2-19n+45}2]$.
Now assume $\Delta(G)=n-4$. Suppose $v_0\in V(G)$, $d(v_0)=n-4$,
$\Gamma(v_0)=\{v_1,\ldots,v_{n-4}\}$ and
$\Gamma_2(v_0)=\{u_1,\ldots,u_t\}$.  Then clearly $t\le n-6$. We may
suppose that $ v_1,\cdots,v_{s_1} $ are all vertices adjacent to
exactly two vertices in $\Gamma_2(v_0)$ and
$v_{s_1+1},\cdots,v_{s_2}$ are all vertices adjacent to exactly one
vertex in $\Gamma_2(v_0)$. Let $V_1=\{v_0,v_1,\ldots,v_{n-4}\}$ and
$V_1'=V(G)-V_1$. By Lemma 4.4, we have $d(v)\le n-5$ for $v\in
V(G)-\{v_0,v_1,\ldots,v_{n-4}\}$.
\par If $s_2\ge n-6$, from (4.1) we see that
 $d(v_i)\le n-5$ for $i=1,2,\ldots,n-6$. Since $d(v)\le n-5$ for all
 $v\in V(G)-\{v_0,\ldots,v_{n-4}\}$ we see that
 $$\align 2e(G)&=\sum_{v\in V(G)}d(v)
 \le d(v_0)+d(v_{n-4})+d(v_{n-5})+(2n-12)(n-5)
 \\&\le 3(n-4)+(2n-12)(n-5)=2n^2-19n+48\endalign $$
 Thus, $e(G)\le [n^2-\f{19}2n+24]$.
 If $s_2<n-6$, then  $s_1\le s_2<n-6$. Using (4.2) we see that
$$\align e(G)&=e(G[V_1])+e(V_1V'_1)+e(G-V_1)
\\&\leq \f
{n^2-7n+12-s_1-2s_2}2+2s_1+s_2-s_1+\b{n-6}2
\\&=\f {2n^2-20n+54+s_1}2
<\f {2n^2-20n+54+n-6}2 =n^2-\f {19}2n+24. \endalign$$ Thus, we
always have $e(G)\le [n^2-\f{19}2n+24]=n^2-10n+24+[\f n2]$.
\par When $n<26$ we have $n^2-\f {19}2n+24<n^2-10n+37=e(K_{n-1}\cup
K_{n-8})$. By the above, $\ex(2n-9;n)=n^2-10n+37$.

\par Now we assume $n\ge 26$. Clearly $e(K_{n-1}\cup K_{n-8})
=n^2-10n+37\le  n^2-\f {19}2n+24$. To prove the result, now we only
need to construct a connected graph $G_0$ of order $2n-9$ such that
$G_0$ does not contain  $T_n^3$ as a subgraph and
$e(G_0)=n^2-10n+24+[\f n2]$. When $n$ is even, we may construct a
regular graph $H$ with degree $n-10$ and
$V(H)=\{v_1,\ldots,v_{n-6}\}$. Let $G_0$ be a graph given by
$V(G_0)=\{v_0,v_1,\ldots,v_{n-4},u_1,\ldots,u_{n-6}\}$ and
$$\align E(G_0)&=E(H)\cup \{v_0v_1,\ldots,v_0v_{n-4},v_1v_{n-5},
\ldots, v_{n-6}v_{n-5}, v_1v_{n-4},\ldots,v_{n-5}v_{n-4},
\\&\q v_1u_1,v_1u_2,v_2u_1,v_2u_2,\ldots,v_{n-7}u_{n-7},v_{n-7}u_{n-6},
v_{n-6}u_{n-7},v_{n-6}u_{n-6},
\\&\q u_1u_2,\ldots, u_1u_{n-6},u_2u_3,\ldots,u_2u_{n-6},
u_3u_{n-6},\ldots,u_{n-7}u_{n-6}\big\}.
\endalign$$
Then $d(v_0)=d(v_{n-5})=d(v_{n-4})=n-4$ and
$d(v_1)=\cdots=d(v_{n-6})=d(u_1)=\cdots=d(u_{n-6})=n-5$. Clearly
$G_0$ does not contain any copies of $T_n^3$ and
$$e(G_0)=\f 12\sum_{v\in V(G_0)}d(v)=
\f{3(n-4)+(2n-12)(n-5)}2=n^2-\f{19}2n+24.$$ When $n$ is odd, let $H$
be a graph with $V(H)=\{v_1,\ldots,v_{n-6}\}$ and
$$E(H)=\{v_1v_2,v_2v_3,\ldots,v_{n-7}v_{n-6},v_{n-6}v_1,v_1v_{\f{n-5}2},
v_2v_{\f{n-3}2},\ldots,v_{\f{n-7}2}v_{n-7}\}.$$
 Then
$d_H(v_1)=\cdots=d_H(v_{n-7})=3$ and $d_H(v_{n-6})=2$.  Let $G_0$ be
a graph with $V(G_0)=\{v_0,v_1,\ldots,v_{n-4},u_1,\ldots,u_{n-6}\}$
and
$$\align &E(G_0)\\&=E(\overline{H})\cup \{v_0v_1,\ldots,v_0v_{n-4},v_1v_{n-5},
\ldots, v_{n-6}v_{n-5}, v_1v_{n-4},\ldots,v_{n-5}v_{n-4},
\\&\q v_1u_1,v_1u_2,v_2u_1,v_2u_2,\ldots,v_{n-8}u_{n-8},v_{n-8}u_{n-7},
v_{n-7}u_{n-8},v_{n-7}u_{n-7}, v_{n-6}u_{n-6},
\\&\q u_1u_2,\ldots, u_1u_{n-6},u_2u_3,\ldots,u_2u_{n-6},
u_3u_{n-6},\ldots,u_{n-7}u_{n-6}\big\}.
\endalign$$
Then  $d(v_0)=d(v_{n-5})=d(v_{n-4})=n-4$,
$d(v_1)=\cdots=d(v_{n-6})=d(u_1)=\cdots=d(u_{n-7})=n-5$ and
$d(u_{n-6})=n-6$. Clearly $G_0$ does not contain any copies of
$T_n^3$ and
$$\align e(G_0)&=\f 12\sum_{v\in V(G_0)}d(v)=\f{3(n-4)+(n-6+n-7)(n-5)+n-6}2
\\&=\f{2n^2-19n+47}2=n^2-10n+24+\big[\f n2\big].\endalign$$
\par By the above, the lemma is proved.

\pro{Theorem 4.4} Let $p,n\in\Bbb N$, $p\ge n\ge 15$ and
$p=k(n-1)+n-8$ with $k\in\Bbb N$. Then
$$ex(p;T_n^3)=\f{(n-2)p-7n+30}2+\max\Big\{\big[\f
n2\big],13\Big\}.$$
\endpro
Proof. By Lemmas 4.1, 2.7 and 4.6,
$$\align \ex(p;T_n^3)&=\f{(n-2)(p-(2n-9))}2+\ex(2n-9;T_n)
\\&=\f{(n-2)(p-(2n-9))}2+n^2-10n+24+\max\Big\{\big[\f
n2\big],13\Big\}
\\&=\f{(n-2)p-7n+30}2+\max\Big\{\big[\f
n2\big],13\Big\}.\endalign$$

\pro{Lemma 4.7} Let $n \in \Bbb N$ with $n\ge 15$. Then
$$ex(2n-8;T_n^3)=n^2-9n+29+\t{max}\Big\{0,\big[\f{n-37}4\big]\Big\}.$$
\endpro
Proof. Let $G\in Ex(2n-8;T_n^3)$. Then clearly $e(G)\ge
e(K_{n-1}\cup K_{n-7})=n^2-9n+29$. Suppose that $G$ is not connected
and $G=G_1\cup \cdots \cup G_s$, where $G_i$ is a component
 of $G$ with $|V(G_i)|=p_i$ and $p_1\le p_2\le \cdots\le p_s$.
By Lemmas 4.1, 2.5 and the argument in the proof of Lemma 2.8, we
have $s=2$ and $p_1\le n-7$. If $p_1=n-7$, then $p_2=n-1$ and so
$e(G)=e(K_{n-7}\cup K_{n-1})=n^2-9n+29$. If $p_1\le n-8$, then
$p_2=2n-8-p_1\ge n$. By Lemmas 4.1 and 2.5, we have $p_1(n-3-p_1)\le
2n-7$. For $3\le p_1\le n-8$ we have $p_1(n-3-p_1)\ge 3(n-6)>2n-7$.
Thus, $p_1=1$ or $2$.   For $p_1=2$ we have $p_2=2n-10=n-1+n-9$. As
$(n-9)(n-8-(n-9))=n-9\ge 6>5+((-1)^n-(-1)^{(n-1)(n-10)})/2$, using
Lemma 4.2 we see that $\Delta(G_2)\le n-5$ and hence
$e(G_2)=\ex(2n-10;K_{1,n-4})=[\f{(2n-10)(n-5)}2]=n^2-10n+25$. Thus,
$e(G)=e(G_1)+2(G_2)=1+n^2-10n+25=n^2-10n+26<n^2-9n+29=e(K_{n-7}\cup
K_{n-1})$. This is impossible. For $p_1=1$ we have $p_2=2n-9$, from
the proof of Lemma 4.6 we see that $n\ge 26$ and
$e(G)=e(G_2)=n^2-10n+24+[\f n2]<n^2-9n+29=e(K_{n-7}\cup K_{n-1})$.
This is also impossible. Therefore, when $G$ is not connected, we
have $e(G)=e(K_{n-7}\cup K_{n-1})=n^2-9n+29.$
\par Assume that $G$ is connected.
By Lemma 4.1, $\Delta(G)\le n-4$. If $\Delta(G)\le n-5$, then
clearly $e(G)=\ex(2n-8;K_{1,n-4})=[\f{(2n-8)(n-5)}2]=n^2-9n+20
<n^2-9n+29=e(K_{n-7}\cup K_{n-1})$. This is a contradiction. Hence
$\Delta(G)=n-4$. Suppose $v_0\in V(G)$, $d(v_0)=n-4$,
$\Gamma(v_0)=\{v_1,\ldots,v_{n-4}\}$ and
$\Gamma_2(v_0)=\{u_1,\ldots,u_t\}$.  Then clearly $t\le n-5$. We may
suppose that $ v_1,\cdots,v_{s_1} $ are all vertices adjacent to
exactly two vertices in $\Gamma_2(v_0)$ and
$v_{s_1+1},\cdots,v_{s_2}$ are all vertices adjacent to exactly one
vertex in $\Gamma_2(v_0)$. Let $V_1=\{v_0,v_1,\ldots,v_{n-4}\}$ and
$V_1'=V(G)-V_1$.   By (4.1), $d(v_i)\le n-5$ for $i=1,2,\ldots,s_2$.
By Lemma 4.4, $d(v)\le n-5$ for $v\in V_1'$. Thus,
$$\align e(G)&=\f 12\sum_{v\in V(G)}d(v)
\le \f{(n-3-s_2)(n-4)+(n-5+s_2)(n-5)}2
\\&=\f{2n^2-17n+37-s_2}2=n^2-9n+29+\f{n-21-s_2}2.\endalign$$
Since $e(G)\ge n^2-9n+29$ we get $s_2\le n-21$.  By (4.2),
$e(G[V_1])\le \f {n^2-7n+12-s_1-2s_2}2$.
  Thus, $$\align
e(G)&=e(G[V_1])+e(V_1V'_1)+e(G-V_1)
\\&\leq \f
{n^2-7n+12-s_1-2s_2}2+2s_1+s_2-s_1+\b {n-5}2 =n^2-9n+21+\f{s_1}2
\\&\le n^2-9n+21+\f{n-21}2=n^2-\f {17}2n+\f{21}2.\endalign$$
As $e(G)\ge n^2-9n+29$ we get $n^2-\f{17}2n+\f{21}2\ge n^2-9n+29$
and so $n\ge 37$.
\par If $s_1\le\f {n-5}2$, from the above we see that
 $$e(G)\le
n^2-9n+21+\f{n-5}4=n^2-9n+29+\f {n-37}4.$$ If $s_1=\f{n-5}2+s_1'>\f
{n-5}2$, then $e(V_1V_1')\ge 2s_1=n-5+2s_1'$. As $d(v)\le n-5$ for
$v\in V_1'$, we see that $e(V_1V_1')+2e(G-V_1)=\sum_{v\in
V_1'}d(v)\le (n-5)(n-5)$ and so $2e(G-V_1)\le (n-5)^2-e(V_1V_1')\le
(n-5)^2-(n-5)-2s_1'=n^2-11n+30-2s_1'$. Therefore,
$$\align e(G)&=e(G[V_1])+e(V_1V'_1)+e(G-V_1)
\\&\le  \f
{n^2-7n+12-s_1-2s_2}2+2s_1+s_2-s_1+\f {n^2-11n+30-2s_1'}2
\\&\q=n^2-9n+29+\f {n-37-2s_1'}4<n^2-9n+29+\f {n-37}4.\endalign$$
Thus, when $G$ is connected, we always have $n\ge 37$ and $e(G)\le
n^2-9n+29+[\f {n-37}4]$.
\par By the above, for $n<37$ we see that $G$ is not connected and
$e(G)=n^2-9n+29$. Now assume $n\ge 37$. Then $n^2-9n+29+[\f
{n-37}4]\ge n^2-9n+29$. To prove the result, we only need to
construct a connected graph $G_0$ of order $n^2-9n+29+[\f {n-37}4]$
without $T_n^3$. Let us consider the following four cases:
\par{\bf Case 1.} $n\e 1\pmod 4$. In this case, by [3, Corollary 2.1]
 we may construct  a regular graph $H$ with degree $\f{n-13}2$ and
$V(H)=\{v_1,\ldots,v_{\f{n-5}2}\}$. Let $G_0$ be a graph with
$V(G_0)=\{v_0,v_1,\ldots,v_{n-4},u_1,\ldots,u_{n-5}\}$ and
$$\align E(G_0)&=E(H)\cup
\big\{u_1u_2,\ldots,u_1u_{n-5},u_2u_3,
\ldots,u_2u_{n-5},u_3u_4,\ldots,u_{n-6}u_{n-5},
\\&\qq
v_1u_1,v_1u_2,\ldots,v_{\f{n-5}2}u_{n-6},v_{\f{n-5}2}u_{n-5},
 v_0v_1,\ldots,v_0v_{n-4}, \\&\qq v_1v_{\f{n-3}2},\ldots,v_1v_{n-4},
\ldots,v_{\f{n-5}2}v_{\f{n-3}2},\ldots,v_{\f{n-5}2}v_{n-4},
\\&\qq v_{\f{n-3}2}v_{\f{n-1}2},\ldots,
v_{\f{n-3}2}v_{n-4},v_{\f{n-1}2}v_{\f{n+1}2},\ldots,v_{n-5}v_{n-4}
\big\}.\endalign$$ Then $d(v_0)=d(v_{\f
{n-3}2})=\cdots=d(v_{n-4})=n-4$ and $d(v_1)=\cdots=d(v_{\f
{n-5}2})=d(u_1)=\cdots=d(u_{n-5})=n-5$. It is clear that $G_0$  does
not contain any copies of $T_n^3$ and
$$\align2e(G_0)&=\sum_{v\in V(G_0)}d(v)=
\Big(n-5+\f{n-5}2\Big)(n-5)+\Big(n-3-\f{n-5}2\Big)(n-4)
\\&=2n^2-18n+58+\f{n-37}2.\endalign$$
Therefore, $e(G_0)=n^2-9n+29+\f {n-37}4$.

\par{\bf Case 2.} $n\e 2\pmod 4$. Let $H$ be  a graph  with $V(H)=\{v_1,\ldots,v_{\f{n-4}2}\}$ and
$$E(H)=\big\{v_1v_2,v_2v_3,\ldots,v_{\f{n-6}2}v_{\f{n-4}2},v_{\f{n-4}2}v_1,
v_1v_{\f{n-2}4},v_2v_{\f{n+2}4},\ldots,
v_{\f{n-6}4}v_{\f{n-6}2}\big\}.$$
 Then
$d_H(v_1)=\cdots=d_H(v_{\f{n-6}2})=3$ and $d_H(v_{\f{n-4}2})=2$. Let
$G_0$ be a graph with
$V(G_0)=\{v_0,v_1,\ldots,v_{n-4},u_1,\ldots,u_{n-5}\}$ and
$$\align E(G_0)&=E(\overline{H})\cup \{u_iu_j\
(i,j\in\{1,2,\ldots,n-5\},i<j),\\&\qq
v_1u_1,v_1u_2,\ldots,v_{\f{n-6}2}u_{n-7},v_{\f{n-6}2}u_{n-6},v_{\f{n-4}2}u_{n-5},
\\&\q\q v_0v_1,\ldots,v_0v_{n-4},v_iv_j
(i\in\{(n-2)/2,\ldots,n-4\},j\in\{1,2,\ldots,i-1\})\}.\endalign$$
Then $d(v_0)=d(v_{\f {n-2}2})=\cdots=d(v_{n-4})=n-4$ and
$d(v_1)=\cdots=d(v_{\f {n-4}2})=d(u_1)=\cdots=d(u_{n-5})=n-5$. It is
clear that $G_0$  does not contain any copies of $T_n^3$ and
$$\align2e(G_0)&=\Big(n-5+\f{n-4}2\Big)(n-5)+\Big(n-3-\f{n-4}2\Big)(n-4)
\\&=2n^2-18n+58+\f{n-38}2.\endalign$$
Therefore, $e(G_0)=n^2-9n+29+[\f {n-37}4]$.

\par{\bf Case 3.} $n\e 3\pmod 4$. Let $H$ be a
graph with $V(H)=\{v_1,\ldots,v_{\f{n-3}2}\}$ and
$$E(H)=\big\{v_1v_2,v_2v_3,\ldots,v_{\f{n-5}2}v_{\f{n-3}2},
v_{\f{n-3}2}v_1, v_1v_{\f{n-3}4},v_2v_{\f{n+1}4},\ldots,
v_{\f{n-7}4}v_{\f{n-7}2}\big\}.$$
 Then
$d_H(v_1)=\cdots=d_H(v_{\f{n-7}2})=3$ and
$d_H(v_{\f{n-5}2})=d_H(v_{\f{n-3}2})=2$.  Let $G_0$ be a graph with
$V(G_0)=\{v_0,v_1,\ldots,v_{n-4},u_1,\ldots,u_{n-5}\}$ and
$$\align E(G_0)&=E(\overline{H})\cup \{u_iu_j\
(i,j\in\{1,2,\ldots,n-5\},i<j),\\&\qq
v_1u_1,v_1u_2,\ldots,v_{\f{n-7}2}u_{n-8},v_{\f{n-7}2}u_{n-7},v_{\f{n-5}2}u_{n-6},v_{\f{n-3}2}u_{n-5},
\\&\q\q v_0v_1,\ldots,v_0v_{n-4},v_iv_j
(i\in\{(n-1)/2,\ldots,n-4\},j\in\{1,2,\ldots,i-1\})\}.\endalign$$
Then $d(v_0)=d(v_{\f {n-1}2})=\cdots=d(v_{n-4})=n-4$ and
$d(v_1)=\cdots=d(v_{\f {n-3}2})=d(u_1)=\cdots=d(u_{n-5})=n-5$.
Clearly $G_0$  does not contain any copies of $T_n^3$ and
$$\align2e(G_0)&=\Big(n-5+\f{n-3}2\Big)(n-5)
+\Big(n-3-\f{n-3}2\Big)(n-4)
\\&=2n^2-18n+58+\f{n-39}2.\endalign$$
Therefore, $e(G_0)=n^2-9n+29+\f{n-39}4=n^2-9n+29+[\f {n-37}4]$.
\par{\bf Case 4.} $n\e 0\pmod 4$. Let $H$ be a
graph with $V(H)=\{v_1,\ldots,v_{\f {n-2}2}\}$ and
$$E(H)=\big\{v_1v_2,v_2v_3,\ldots,v_{\f {n-4}2}v_{\f {n-2}2},
v_{\f {n-2}2}v_1, v_1v_{\f {n-4}4},v_2v_{\f n4},\ldots,v_{\f
{n-8}4}v_{\f {n-8}2}\big\}.$$
 Then
$d_H(v_1)=\cdots=d_H(v_{\f{n-8}2})=3$ and
$d_H(v_{\f{n-6}2})=d_H(v_{\f{n-4}2})=d_H(v_{\f{n-2}2})=2$. Let $G_0$
be a graph with
$V(G_0)=\{v_0,v_1,\ldots,v_{n-4},u_1,\ldots,u_{n-5}\}$ and
$$\align E(G_0)&=E(\overline{H})\cup \{u_iu_j\
(i,j\in\{1,2,\ldots,n-5\},i<j),\\&\qq
v_1u_1,v_1u_2,\ldots,v_{\f{n-8}2}u_{n-9},v_{\f{n-8}2}u_{n-8},v_{\f{n-6}2}u_{n-7},v_{\f{n-4}2}u_{n-6},v_{\f{n-2}2}u_{n-5},
\\&\q\q v_0v_1,\ldots,v_0v_{n-4},v_iv_j
(i\in\{n/2,\ldots,n-4\},j\in\{1,2,\ldots,i-1\})\}.\endalign$$ Then
$d(v_0)=d(v_{\f n2})=\cdots=d(v_{n-4})=n-4$ and
$d(v_1)=\cdots=d(v_{\f {n-2}2})=d(u_1)=\cdots=d(u_{n-5})=n-5$.
Clearly $G_0$  does not contain any copies of $T_n^3$ and
$$\align2e(G_0)&=\Big(n-5+\f{n-2}2\Big)(n-5)+\Big(n-3-\f{n-2}2\Big)(n-4)
\\&=2n^2-18n+58+\f{n-40}2.\endalign$$
Therefore, $e(G_0)=n^2-9n+29+\f {n-40}4=n^2-9n+29+[\f{n-37}4]$.
\par Summarizing the above we prove the lemma.

\pro{Theorem 4.5} Let $p,n\in\Bbb N$, $p\ge n\ge 15$ and
$p=k(n-1)+n-7$ with $k\in\Bbb N$. Then
$$ex(p;T_n^3)=\f{(n-2)p-6(n-7)}2+\max\Big\{\big[\f
{n-37}4\big],0\Big\}.$$
\endpro
Proof.  By Lemmas 4.1, 2.7 and 4.7,
$$\align \ex(p;T_n^3)&=\f{(n-2)(p-(2n-8))}2+\ex(2n-8;T_n)
\\&=\f{(n-2)(p-(2n-8))}2+n^2-9n+29+\max\Big\{\big[\f
{n-37}4\big],0\Big\}
\\&=\f{(n-2)p-6(n-7)}2+\max\Big\{\big[\f
{n-37}4\big],0\Big\}.\endalign$$

\section*{5. Evaluation of $ex(p;T_n''')$}
\pro{Lemma 5.1} Let $p,n \in \Bbb N,\ p\geq n\geq 10$ and $G\in
Ex(p;T_n''')$. Suppose  that G is connected. Then $\Delta(G)\le
n-4$.\endpro Proof.  By Lemma 2.1, $\Delta(G)\geq n-5$. Suppose
$v_0\in V(G),d(v_0)=\Delta(G)=m$ and
$\Gamma(v_0)=\{v_1,\ldots,v_m\}$. If $m=p-1$, as G does not contain
$T_n'''$ as a subgraph, we see that $G[v_1,\ldots,v_m]$ does not
contain $3K_2$ as a subgraph. If $G[v_1,\ldots,v_m]$ does not
contain $2K_2$ as a subgraph, then G does not contain $T_n'''$ as a
subgraph, $e(G[v_1,\ldots,v_m])\le e(K_{1,m-1})=m-1$ and so $e(G)=
d(v_0)+e(G[v_1,\ldots,v_m])\le m+m-1=2m-1$. Suppose that
$G[v_1,\ldots,v_m]$ contains $2K_2$ as a subgraph and
$v_1v_2,v_3v_4\in E(G)$. Then every edge in $G[v_1,\ldots,v_m]$ is
incident with some vertex in $\{v_1,v_2,v_3,v_4\}$. If
$v_2v_i,v_3v_i,v_4v_i\in E(G)$ for some $i\in\{5,\ldots,m\}$, then
all edges in $E(G[v_1,\ldots,v_m])-\{v_3v_4\}$ is incident with
$v_2$ or $v_i$ and so $e(G)=d(v_0)+e(G[v_1,\ldots,v_m])\le
m+d(v_2)-1+d(v_i)-1\le 3m-2$. If $d(v_i)\le 3$ for $i=5,\ldots,m$,
then there are at most two vertices in $\{v_1,v_2,v_3,v_4\}$
adjacent to a fixed vertex in $\{v_5,\ldots,v_m\}$ and so $e(G)\le
d(v_0)+2(m-4)+2=3m-6.$ From the above we always have $e(G)\le 3m$.
This contradicts Lemma 2.2. Hence $m<p-1$. Suppose that
$u_1,\ldots,u_t$ are all vertices in G such that
$d(u_1,v_0)=\cdots=d(u_t,v_0)=2$. Then $t\geq1$. We may assume
$u_1v_1\in E(G)$ with no loss of generality. Set
$V_1=\{v_0,v_1,\ldots,v_m\}$ and $V_2=\{v_0,v_1,\ldots,v_m,u_1\}$.
\par Suppose $t=1$ and $m\ge n-2$. Let $G'=G[v_2,v_3,\ldots,v_m]$.
 As $G$ does not contain
$T_n'''$, we see that $G'$ does not contain any copies of $2K_2$. If
$e(G')\le 2$, then
$$e(G)-e(G-V_1)\le d(v_0)+d(u_1)+d(v_1)-2+e(G')\le m+m+m-2+2=3m,$$
which contradicts Lemma 2.2. Hence $e(G')\ge 3$. We claim that $G'$
does not contain any copies of $K_3$. We may assume $v_2v_3,v_2v_4,
v_3v_4\in E(G')$. As $G'$ does not contain any copies of $2K_2$ we
see that $e(G')=3$. If $|\Gamma(u_1)\cap\{v_2,\ldots,v_m\}|\ge 2$
and $|\Gamma(u_1)\cap\{v_5,\ldots,v_m\}|\ge 1$,  then $v_1$ is not
adjacent to any vertex in $G'$. Hence $e(G)-e(G-V_1)\le
d(v_0)+d(u_1)+e(G')\le m+m+3<3m$. By Lemma 2.2, this is impossible.
Similarly, if $|\Gamma(v_1)\cap\{v_2,\ldots,v_m\}|\ge 2$ and
$|\Gamma(v_1)\cap\{v_5,\ldots,v_m\}|\ge 1$, then $u_1$ is not
adjacent to any vertex in $G'$. Hence $e(G)-e(G-V_1)\le
d(v_0)+d(v_1)-1+e(G')\le m+m-1+3<3m$, which contradicts Lemma 2.2.
As $m\ge n-2\ge 8$, $d(v_1)\le m$ and $d(u_1)\le m$, from the above
we
 have
$$|\Gamma(v_1)\cap V(G')|+|\Gamma(u_1)\cap V(G')|\le \max
\{3+3,m-1\}=m-1.$$ Thus, $$e(G)-e(G-V_1)\le
d(v_0)+1+|\Gamma(v_1)\cap V(G')|+|\Gamma(u_1)\cap V(G')|\le
m+1+m-1=2m.$$ This is also impossible by Lemma 2.2. Hence the claim
is true.
\par Now assume that $e(G')\ge 3$ and $G'$ does not contain any
copies of $2K_2$ and $K_3$. Then all edges in $G'$ have a common
endpoint. We may assume that $v_2$ is  such a vertex.
 Therefore $d_{G'}(v_2)\ge 3$. Suppose $v_1v_i\in E(G)$
for some $i\in\{3,4,\ldots,m\}$. Then $u_1v_j\not\in E(G)$ for all
$j\in\{3,4,\ldots,m\}-\{i\}$ and so $|\Gamma(u_1)\cap
\{v_2,\ldots,v_m\}|\le 2$. Otherwise, for some $v_k\in \Gamma(v_2)$
the three edges $v_1v_i,u_1v_j,v_2v_k$ induce a copy of $3K_2$ and
so $G$ contains a copy of $T_n'''$. Hence $d(v_1)+|\Gamma(u_1)\cap
\{v_2,\ldots,v_m\}|\le m+2$. If $v_1v_i\not\in E(G)$ for every
$i\in\{3,4,\ldots,m\}$, then $d(v_1)\le 3$ and
$d(v_1)+|\Gamma(u_1)\cap \{v_2,\ldots,v_m\}|\le 3+(m-1)=m+2$.
Therefore,
$$\align e(G)-e(G-V_1)&=d(v_0)-1+d(v_1)+|\Gamma(u_1)\cap
\{v_2,\ldots,v_m\}|+d_{G'}(v_2) \\&\le
m-1+(m+2)+(m-2)<3m.\endalign$$ This is impossible by Lemma 2.2.
Hence $\Delta(G)\le n-3$ for $t=1$.
\par Suppose $t=1$ and $\Delta(G)=m\in\{n-3,n-4\}$. Then
$$\align e(G)-e(G-V_2)&=d(u_1)+e(G[v_0,v_1,\ldots,v_m])
\\&\le m+e(K_{m+1})=\f{m^2+3m}2<\f{(m+1)(m+2)}2=e(K_{m+2}).
\endalign$$
Thus, $e(G)<e((G-V_2)\cup K_{m+2})$, which contradicts the
assumption $G\in Ex(p;T_n''')$.

By the above, for $t=1$ we have $\Delta(G)\le n-5$. From now on we
assume that $t\ge 2$. Suppose $t=2$, $u_1v_1,u_2v_2\in E(G)$ and
$m=\Delta(G)\ge n-3$. As $G$ does not contain any copies of
$T_n'''$, we see that $\{v_3,\ldots,v_m\}$ is an independent set in
$G'$. If $i,j\in\{3,4,\ldots,m\}$, $i\not=j$ and $v_1v_i,u_1v_j\in
E(G)$, then $u_2v_2,v_1v_i,u_1v_j$ induce a copy of $3K_2$ and so
$G$ contains a copy of $T_n'''$. Hence
$|\Gamma(v_1)\cap\{v_3,\ldots,v_m\}|
+|\Gamma(u_1)\cap\{v_3,\ldots,v_m\}|\le m-2$. Similarly,
$|\Gamma(v_2)\cap\{v_3,\ldots,v_m\}|
+|\Gamma(u_2)\cap\{v_3,\ldots,v_m\}|\le m-2$. If $u_1v_r,u_2v_s\in
E(G)$, where $r,s\in\{3,4,\ldots,m\}$ and $r\not=s$, then
$v_1v_2\not\in E(G)$, otherwise $u_1v_r,u_2v_s,v_1v_2$ induce a copy
of $3K_2$ and so $G$ contains a copy of $T_n'''$. Hence

$$\align e(G)-e(G-V_1)&=d(v_0)
+e(G[v_1,v_2,u_1,u_2])-e(G([u_1,u_2])
\\&\q +|\Gamma(v_1)\cap\{v_3,\ldots,v_m\}|
+|\Gamma(u_1)\cap\{v_3,\ldots,v_m\}|
\\&\q+|\Gamma(v_2)\cap\{v_3,\ldots,v_m\}|
+|\Gamma(u_2)\cap\{v_3,\ldots,v_m\}|
\\&\le m+4+(m-2)+(m-2)=3m.\endalign$$
By Lemma 2.2, $e(G)-e(G-V_1)>3m$. We get a contradiction. Hence
$\Delta(G)=m\le n-4$.
\par Suppose $t=2$ and $v_1u_1,v_1u_2\in E(G)$. If $u_iv_j\in
E(G)$ for some $i\in\{1,2\}$ and $j\in\{2,3,\ldots,m\}$, by the
above argument we have $\Delta(G)\le n-4$. Now suppose that $u_iv_j
\notin E(G)$ for every $i=1,2$ and $j=2,3,\ldots,m$.
 If $m\ge n-2$, then $G'$ does not contain $2K_2$ as a
subgraph. If $G'$ contains a copy of $K_3$, then $e(G')=3$ and so
$$e(G)-e(G-V_1)=d(v_0)+d(v_1)-1+e(G')\le m+m-1+3=2m+2.$$ This is
impossible by Lemma 2.2. Thus all edges in $G'$ have a common
endpoint and so $e(G')\le e(K_{1,m-2})=m-2$. Hence
 $e(G)-e(G-V_1)=d(v_0)+d(v_1)-1+e(G')\le
m+m-1+m-2=3m-3$. By Lemma 2.2, this is impossible. Therefore
$m=\Delta(G)\le n-3$. If $m=n-3$, then
$$\align
e(G)-e(G-V_2)&=d(v_1)+d(u_1)-1+e(G[v_0,v_2,v_3,\ldots,v_{n-3}])
\\&\le n-3+n-3-1+\binom {n-3}2=\f{n^2-3n-2}2
\\&<\f {(n-1)(n-2)}2=e(K_{n-1}).\endalign$$
Thus, $e(G)<e((G-V_2)\cup K_{n-1})$, which contradicts the
assumption $G\in Ex(p;T_n''')$. Therefore $\Delta(G)\le n-4$.
\par From now on we assume $t\ge 3$.
Suppose  $|\Gamma(v_1)\cap \Gamma_2(v_0)|\ge 2$ and
$|\Gamma(v_2)\cap \Gamma_2(v_0)|\ge 1$. If $|\Gamma(v_2)\cap
\Gamma_2(v_0)|=1$ and $v_2u_2\in E(G)$, then $\{v_3,\ldots,v_m\}$ is
an independent set in $G'$ and $u_iv_j\notin E(G)$ for any $i\in
\{1,3,\ldots,t\}$ and $j\in\{3,4,\ldots,m\}$. Suppose $v_2v_i\in
E(G)$ for some $i\in\{3,4,\ldots,m\}$. Then $u_2v_j\not\in E(G)$ for
all $j\in\{3,4,\ldots,m\}-\{i\}$ and so $|\Gamma(u_2)\cap
\{v_2,\ldots,v_m\}|\le 2$. Otherwise, the three edges
$v_2v_i,u_2v_j,u_1v_1$ induce a copy of $3K_2$ and so $G$ contains a
copy of $T_n'''$. Hence $d(v_2)+|\Gamma(u_2)\cap
\{v_2,\ldots,v_m\}|\le m+2$. If $v_2v_i\not\in E(G)$ for every
$i\in\{3,4,\ldots,m\}$, then $d(v_2)\le 3$ and
$d(v_2)+|\Gamma(u_2)\cap \{v_2,\ldots,v_m\}|\le 3+(m-1)=m+2$.
 Hence
$$e(G)-e(G-V_1)=d(v_0)+d(v_1)+d(v_2)-2+|\Gamma(u_2)\cap \{v_2,\ldots,v_m\}|\le 3m.$$
If $|\Gamma(v_1)\cap\Gamma_2(v_0)|\ge 2$ and $|\Gamma(v_2)\cap
\Gamma_2(v_0)|\ge 2$, then $\{v_3,\ldots,v_m\}$ is an independent
set in $G'$ and $u_iv_j\notin E(G)$ for any $i\in \{1,2,\ldots,t\}$
and $j\in\{3,4,\ldots,m\}$. Hence
$$e(G)-e(G-V_1)=d(v_0)+d(v_1)+d(v_2)-2\le 3m-2.$$ From the above we
always have $e(G)-e(G-V_1)\le 3m$, which contradicts Lemma 2.2.
Therefore $m=\Delta(G)\le n-4$.

\pro{Lemma 5.2} Let $p,n \in \Bbb N,\ p\geq n\geq 10$ and $G\in
Ex(p;T_n''')$. Suppose  that G is connected. Then
$\Delta(G)=n-5$.\endpro Proof. By Lemma 5.1, $\Delta(G)\le n-4$.
 Suppose that $v_0\in V(G), d(v_0)=\Delta(G)=n-4$,
$\Gamma(v_0)=\{v_1,\ldots,v_{n-4}\}$ and
$\Gamma_2(v_0)=\{u_1,\ldots,u_t\}$. By the proof of Lemma 5.1, we
have $\Delta(G)=n-5$ for $t=1$. Set
$V_1=\{v_0,v_1,\ldots,v_{n-4},u_1,u_2\}$.

Suppose $t=2$. Then
$$\align e(G)-e(G-V_1)&=d(u_1)+d(u_2)+e(G[v_0,v_1,\ldots,v_{n-4}])
\\&\le n-4+n-4+\f{(n-3)(n-4)}2\\&=\f{n^2-3n-4}2<\f{n^2-3n+2}2=e(K_{n-1}).\endalign$$
Thus, $e(G)<e((G-V_1)\cup K_{n-1})$, which contradicts the
assumption $G\in Ex(p;T_n''')$.

Now assume $t\ge 3$. If $|\Gamma(v_1)\cap \Gamma_2(v_0)|=t$, then
$u_iv_j\notin E(G)$ for any $i\in \{1,2,\ldots,t\}$ and
$j\in\{2,4,\ldots,n-4\}$. We see that
$$\align e(G)-e(G-V_1)&\le
d(v_1)+d(u_1)+d(u_2)+e(G[v_0,v_2,\ldots,v_{n-4}])
\\&\le n-4+n-4+n-4+\f{(n-4)(n-5)}2
=\f{n^2-3n-4}2\\&<\f{(n-1)(n-2)}2=e(K_{n-1})\endalign$$ and so
$e(G)<e((G-V_1)\cup K_{n-1})$. This contradicts the assumption $G\in
Ex(p;T_n''')$. If $|\Gamma(v_1)\cap \Gamma_2(v_0)|\ge 1$ and
$|\Gamma(v_2)\cap \Gamma_2(v_0)|\ge 1$, then $u_1v_1,u_2v_2\in
E(G)$, $u_iv_j\notin E(G)$ for any $i\in \{3,4,\ldots,t\}$ and
$j\in\{3,4,\ldots,n-4\}$. Thus,  $$\align e(G)-e(G-V_1)&\le
d(v_1)+d(v_2)+d(u_1)+d(u_2)+e(G[v_0,v_3,\ldots,v_{n-4}])
\\&\le n-4+n-4+n-4++n-4+\f{(n-5)(n-6)}2
\\&=\f {n^2-3n-2}2<\f {(n-1)(n-2)}2=e(K_{n-1}).\endalign$$
This contradicts the fact $G\in Ex(p;T_n''')$. Hence $\Delta(G)=n-5$
as claimed.

\pro{Theorem 5.1} Let $p,n \in \Bbb N, p\geq n\geq 10$,
$p=k(n-1)+r$, $k\in\Bbb N$ and $r\in\{0,1,\ldots,n-2\}$. Then
$$ex(p;T_n''')=\f{(n-2)p-r(n-1-r)}2+\max\Big\{0,\big[\f{r(n-4-r)-3(n-1)}2\big]\Big\}.$$\endpro
Proof. This is immediate from Lemmas 2.10 and 5.2.

\end{CJK*}

\begin{thebibliography}{99}

 \bibitem{1} R.J. Faudree and R.H. Schelp,
Path Ramsey numbers in multicolorings,  J. Combin. Theory Ser. B
{\bf 19}(1975), 150-160.


 \bibitem{2} A.F. Sidorenko, Asymptotic solution for a new
class of forbidden $r$-graphs,  Combinatorica {\bf  9}(1989),
207-215.

 \bibitem{3} Z.H. Sun and L.L.Wang,
Tur\'an's problem for trees, J. Combin. Number Theory {\bf 3}(2011),
51-69.
 \bibitem{4} Z.H. Sun, L.L.Wang and Y.L. Wu, Tur\'an's problem and Ramsey
 numbers for trees, arXiv:1110.2725.
 \bibitem{5} M. Wo\'zniak, On the Erd\H os-S\'os
conjecture, J. Graph Theory {\bf 21}(1996), 229-234.
\end{thebibliography}
\end{document}